\documentclass[preprint,12pt]{elsarticle}

\usepackage{amssymb}
\usepackage{amsmath}
\usepackage{amsthm}
\usepackage{graphicx} 
\usepackage[font=small,labelfont=bf]{caption}
\usepackage{subcaption}
\usepackage{multirow}
\usepackage{booktabs}
\usepackage{algorithm}
\usepackage{algpseudocode}
\usepackage{hyperref} 
\usepackage{array}
\usepackage{caption}

\newtheorem{proposition}{Proposition}

\newtheorem{remark}{Remark}

\journal{Computational and Applied Mathematics}
\begin{document}
\begin{frontmatter}



\title{A New Combination of Preconditioned Gradient Descent Methods and Vector Extrapolation Techniques for Nonlinear Least-Squares Problems}


\author[inst1,inst2]{Abdellatif Mouhssine}
\affiliation[inst1]{organization={The UM6P Vanguard Center, Mohammed VI Polytechnic University},
             addressline={Lot 660 Hay Moulay Rachid}, 
            city={Benguerir},
            postcode={43150}, 
            country={Morocco}}
\ead{abdellatif.mouhssine@um6p.ma}           


\affiliation[inst2]{organization={Laboratory of Pure and Applied Mathematics, University of Littoral C\^ote d'Opale},
             addressline={50 Rue F. Buisson}, 
            city={Calais Cedex},
            postcode={B.P. 699 - 62228}, 
            country={France}}

\begin{abstract}
Vector extrapolation methods are widely used in large-scale simulation studies, and numerous extrapolation-based acceleration techniques have been developed to enhance the convergence of linear and nonlinear fixed-point iterative methods. While classical extrapolation strategies often reduce the number of iterations or the computational cost, they do not necessarily lead to a significant improvement in the accuracy of the computed approximations. In this paper, we study the combination of preconditioned gradient-based methods with extrapolation strategies and propose an extrapolation-accelerated framework that simultaneously improves convergence and approximation accuracy. The focus is on the solution of nonlinear least-squares problems through the integration of vector extrapolation techniques with preconditioned gradient descent methods. A comprehensive set of numerical experiments is carried out to study the behavior of polynomial-type extrapolation methods and the vector $\varepsilon$-algorithm when coupled with gradient descent schemes, with and without preconditioning. The results demonstrate the impact of extrapolation techniques on both convergence rate and solution accuracy, and report iteration counts, computational times, and relative reconstruction errors. The performance of the proposed hybrid approaches is further assessed through a benchmarking study against Gauss--Newton methods based on generalized Krylov subspaces.

\end{abstract}




\end{frontmatter}


\section{Introduction}\label{sec1}

Extrapolation methods \cite{smith1987extrapolation,eddy1979extrapolating,mesina1977convergence,cabay1976apolynomial,sidi1986convergence,brezinski1991extrapolation,brezinski2000convergence} and acceleration techniques \cite{saad2025acceleration} are frequently employed because the convergence of sequences generated by linear or nonlinear fixed-point iterative methods is often too slow. Vector extrapolation methods were specifically designed to handle such vector sequences and have been successfully applied in a wide range of scientific computing problems. These methods include vector polynomial extrapolation techniques \cite{jbilou2000vector,jbilou1995analysis,jbilou1991some}, vector and topological epsilon algorithms \cite{brezinski1991extrapolation,brezinski2000convergence}, and Anderson acceleration methods \cite{anderson1965iterative,walker2011anderson,walker2011andersonalgo,fang2009two,toth2015convergence}.

Among vector polynomial extrapolation methods, the most well-known approaches are Aitken’s method \cite{aitken1927xxv}; the reduced rank extrapolation (RRE) method, attributed to Eddy \cite{eddy1979extrapolating} and Mesina \cite{mesina1977convergence}; the minimal polynomial extrapolation (MPE) method, introduced by Cabay and Jackson \cite{cabay1976apolynomial}; and the modified minimal polynomial extrapolation (MMPE) method, proposed by Sidi, Ford, and Smith \cite{sidi1986acceleration}, Brezinski \cite{brezinski1975generalisations}, and Pugatchev \cite{pugachev1977acceleration}. The epsilon algorithms (or $\epsilon$-algorithms) include the scalar and vector epsilon algorithms (SEA and VEA) of Wynn \cite{wynn1956device,wynn1962acceleration}, as well as the topological epsilon algorithm (TEA) of Brezinski \cite{brezinski1975generalisations}. These topological algorithms have recently been simplified and are now commonly referred to as simplified topological $\epsilon$-algorithms \cite{brezinski2014simplified}. Several convergence results and theoretical properties of the different classes of extrapolation methods can be found in \cite{brezinski1975generalisations,jbilou1995analysis,sadok1993quasilinear,sidi1986convergence,sidi1986acceleration,sidi1988extrapolation,brezinski2018shanks,brezinski2022shanks}.

For linearly generated sequences, polynomial extrapolation techniques are closely related to Krylov subspace methods \cite{sidi1988extrapolation}. In particular, the RRE method can be interpreted as an orthogonal projection method and is equivalent to GMRES \cite{saad1981Krylov}, while the MPE method corresponds to an oblique projection method equivalent to the Arnoldi algorithm \cite{saad1981Krylov}. As in Krylov subspace methods, the computational cost and storage requirements of extrapolation techniques increase with the number of fixed-point iterates involved. Consequently, restarting strategies \cite{mouhssine2025vector,mouhssine2025enhancing} are often employed to control both memory usage and computational effort. This leads to what are commonly referred to as cyclic or restarted extrapolation approaches. In this work, we consider two practical strategies for applying extrapolation algorithms: the Acceleration Method (AM) and the Restarted Method (RM).

In general, vector extrapolation methods are effective for both linear and nonlinear problems when the underlying fixed-point iterates exhibit a sufficiently regular error structure. Such behavior typically arises in problems with asymptotically linear convergence, smooth error decay, or elliptic-type models; see \cite{mouhssine2025vector,mouhssine2025enhancing}. Similar characteristics are frequently observed in fixed-point iterations for nonlinear equations.

Building on our recent works~\cite{mouhssine2025vector,mouhssine2025enhancing}, where vector extrapolation methods were shown to be effective for accelerating linear and nonlinear fixed-point iterations, the present study extends these approaches to a broader and practically relevant class of problems, namely nonlinear least-squares problems, with a particular focus on gradient-based fixed-point iterations and the role of problem conditioning.

The solution of nonlinear least-squares problems \cite{buccini2023efficient} is addressed here by combining vector extrapolation techniques with preconditioned gradient-based methods. In particular, polynomial-type extrapolation strategies (RRE and MPE) and the vector epsilon algorithm (VEA) are integrated with preconditioned gradient descent schemes, including PGD and SGD. A comprehensive numerical study is carried out on several test problems, including the classical Bratu problem \cite{keller1987lectures} in its ill-conditioned regime (large values of the parameter $\lambda$), an extended formulation of the Bratu problem, and an extremely sparse nonlinear problem. Compared to classical gradient-based solvers, the proposed hybrid approaches lead to a significant reduction in the relative error, a decrease in the number of iterations, and an overall acceleration of convergence. In particular, the numerical results show that the proposed methods outperform Gauss–-Newton approaches based on generalized Krylov subspaces \cite{buccini2023efficient}; see also \cite{lampe2012large,huang2017majorization,lanza2015generalized,buccini2022fast,buccini2023limited} for definitions and applications of generalized Krylov subspaces. These results highlight the potential of vector extrapolation techniques as a competitive and effective tool for solving large-scale nonlinear problems.

The remainder of this paper is organized as follows. In Section \ref{sec2}, we introduce vector extrapolation methods, including polynomial-type extrapolation techniques such as reduced rank extrapolation (RRE) and minimal polynomial extrapolation (MPE), as well as the vector $\varepsilon$-algorithm (VEA), and discuss their implementation in detail. Section \ref{sec3} reviews iterative methods for nonlinear least-squares problems and presents an accelerated approach based on the application of polynomial-type extrapolation techniques and the vector $\varepsilon$-algorithm to preconditioned gradient descent schemes, highlighting the role of preconditioning in improving both convergence and robustness. Section \ref{sec4} reports a series of numerical experiments, including results for the RRE-PGD, MPE-PGD, RRE-SGD, and MPE-SGD methods, together with comparisons against the VEA-PGD, VEA-SGD, and Gauss–-Newton methods based on generalized Krylov subspaces.

\section{Vector extrapolation methods}   \label{sec2}
In this section, we will first review various vector extrapolation methods \cite{mouhssine2025vector,mouhssine2025enhancing} to apply them to preconditioned gradient descent methods.

\subsection{Polynomial extrapolation methods} \label{subsec:descriptpolyextrapol}
Let $(s_k)_{k\in \mathbb{N}}$ be a sequence of vectors of $\mathbb{R}^{N}$. Polynomial extrapolation methods, namely RRE and the MPE, when applied to the sequence $(s_k)$, produce an approximation $t_{k,q}$ of the limit or the anti-limit of $(s_k)$, where $k$ is the index of the first vector in the sequence, and $q$ is the number of terms used in the extrapolation. This approximation is defined by
\begin{equation}\label{tkq}
t_{k,q}= \sum_{j=0}^{q} \gamma_{j} s_{k+j},
\end{equation}
where 
\begin{equation}\label{gammaj}
    \sum_{j=0}^{q} \gamma_{j}=1 \;\; \text{and}\; \sum_{j=0}^{q}\eta_{i,j}^{(k)} \gamma_{j}=0,\; i=0,\ldots,q-1,\end{equation} 
with the scalars $\eta_{i,j}^{(k)}$ defined by $\eta_{i,j}^{(k)}= (y_i^{(k)}, \Delta s_{k+j})$, where $(\cdot, \cdot)$ denotes the Euclidean inner product and $y_i^{(k)} \in \mathbb{R}^{N}$ defines the method used.\\
$\Delta s_{k}$ and $\Delta^{2}s_{k}$ denote, respectively, the first and second forward differences of $s_k$ and are defined by

\begin{equation*}
\begin{array}{rcl}
\Delta s_{k} & = & s_{k+1}- s{_k},\;\;\;\;\;\;\; \;\;k=0,1,\ldots\\
\Delta ^ 2 s_{k}& = &  \Delta s_{k+1}-\Delta s_{k} ,\;\;\;k=0,1,\ldots
\end{array}
\end{equation*}
The transformation (\ref{tkq}) can be represented as a quotient of two determinants by using (\ref{gammaj})
\begin{equation}\label{quotientdeterminants}
t_{k,q} =
\frac{
\begin{vmatrix}
s_k & s_{k+1} & \cdots & s_{k+q} \\
\eta_{0,0} & \eta_{0,1} & \cdots & \eta_{0,q} \\
\vdots & \vdots & \ddots & \vdots \\
\eta_{q-1,0} & \eta_{q-1,1} & \cdots & \eta_{q-1,q}
\end{vmatrix}
}{
\begin{vmatrix}
1 & 1 & \cdots & 1 \\
\eta_{0,0} & \eta_{0,1} & \cdots & \eta_{0,q} \\
\vdots & \vdots & \ddots & \vdots \\
\eta_{q-1,0} & \eta_{q-1,1} & \cdots & \eta_{q-1,q}
\end{vmatrix}
}.
\end{equation}
Let $Y_q$ be the matrix formed by the vectors $y_0^{(k)},\cdots , y_{q-1}^{(k)}$ and let us introduce the matrices
$$ \Delta ^{i}S_{k,q}=[\Delta ^{i}s_{k},\cdots,\Delta ^{i}s_{k+q-1}], \,i=1,2.$$
If we replace, in the numerator and denominator of (\ref{quotientdeterminants}), each column $j, j = q +1,q,\cdots ,2$ by the
difference with column $j-1$, we obtain the following expression

\begin{equation}\label{newquotientdeterminants}
t_{k,q} =
\frac{
\begin{vmatrix}
s_k & \Delta S_{k,q}  \\
Y_{q}^{T} \Delta s_{k} & Y_{q}^{T} \Delta^{2} S_{k,q} 
\end{vmatrix}
}{
\begin{vmatrix}
Y_{q}^{T} \Delta^{2} S_{k,q}
\end{vmatrix}
}.
\end{equation}
Using Schur’s formula, the approximation $t_{k,q}$ can be written in matrix form as 
\begin{equation}
t_{k,q}= s_{k}-\Delta S_{k,q}\left({ Y_{q}}^{T} \Delta^{2} S_{k,q} \right)^{-1} { Y_{q}}^{T}\Delta s_{k}.
    \label{eq:extrapolapprox}
\end{equation}
\begin{proposition}
Let $t_{k,q}$ be the transformation given by formula (\ref{eq:extrapolapprox}). Then, the transformation exists and is unique if and only if $\mathrm{det}\left({ Y_{q}}^{T} \Delta^{2} S_{k,q} \right) \not = 0$.
\end{proposition}
\noindent We define a new approximation $\tilde{t}_{k,q}= \sum_{j=0}^{q} \gamma_{j} s_{k+j+1},$ of the limit or the anti-limit of $s_k, k\rightarrow \infty$. The generalized residual of $t_{k,q}$ has been defined as 

\begin{equation}
\begin{array}{rcl}
\tilde{r}(t_{k,q}) & = & \tilde{t}_{k,q}- t_{k,q}. \\
\end{array}
\label{eq:resgeneralized}
\end{equation} 
Using (\ref{eq:extrapolapprox}), this residual can be expressed as
\begin{equation}
\begin{array}{rcl}
\tilde{r}(t_{k,q}) & = &  \Delta s_{k}-\Delta^{2} S_{k,q}\left({ Y_{q}}^{T} \Delta^{2} S_{k,q} \right)^{-1} { Y_{q}}^{T}\Delta s_{k}.
\end{array}
\label{eq:newresgeneralized}
\end{equation} 

\begin{proposition}
If we define the matrices $ W_{k,q}=\mathrm{span}\{\Delta^{2} s_{k},\cdots ,\Delta ^{2}s_{k+q-1}\}$ and $L_{k,q} =\mathrm{span}\{y_0^{(k)},\cdots,y_{q-1}^{(k)}\}$, according to (\ref{eq:newresgeneralized}), we show that the generalized residual satisfies
 $$\begin{cases}
\tilde{r}(t_{k,q}) - \Delta s_{k} \in W_{k,q}, \\
\tilde{r}(t_{k,q}) \perp L_{k,q}.
  \end{cases}$$
\end{proposition}

\begin{remark}
These last two conditions show that the generalized residual is obtained by
projecting the vector $\Delta s_k$ onto the subspace $W_{k,q}$ orthogonally to $L_{k,q}$.   
\end{remark}

\begin{remark}
If we assume that the sequence $(s_k)_{k\in \mathbb{N}}$ is produced by a linear or nonlinear fixed-point iteration, starting with a vector $s_0$,
\begin{equation}\label{fixedpoint}
s_{k+1}= G(s_{k}),\; k=0,1,2,\ldots
\end{equation}
whose solution is denoted by $s$, and where $G : \mathbb{R}^N \to \mathbb{R}^N$ is the fixed-point operator. It is mentioned in \cite{sidi1991efficient} that the vector $\tilde{r}(t_{k,q})$ is
\begin{itemize}
    \item  The exact residual for
linear systems $s = G(s)$.
    \item A first-order approximation to the exact residual for nonlinear systems $s = G(s)$.
\end{itemize}
\end{remark}

\noindent We will assume in the following that $k$ is fixed. Without any limitations, we will even assume that $k = 0$ and we will
replace the notation $t_{0,q}$ with $t_q$. The linear system (\ref{gammaj}) is expressed as 
\begin{equation}\label{eqgamma}
\left\{
\begin{array}{ccccccccc}
\gamma_0 &+& \gamma_1 &+& \cdots &+& \gamma_q &=& 1 \\[6pt]
\gamma_0 (y_0,\Delta s_0) &+& \gamma_1 (y_0,\Delta s_1) &+& \cdots &+& \gamma_q (y_0,\Delta s_q) &=& 0 \\[6pt]
\gamma_0 (y_1,\Delta s_0) &+& \gamma_1 (y_1,\Delta s_1) &+& \cdots &+& \gamma_q (y_1,\Delta s_q) &=& 0 \\[6pt]
\vdots & & \vdots & & \ddots & & \vdots & & \vdots \\[6pt]
\gamma_0 (y_{q-1},\Delta s_0) &+& \gamma_1 (y_{q-1},\Delta s_1) &+& \cdots &+& \gamma_q (y_{q-1},\Delta s_q) &=& 0.
\end{array}
\right.
\end{equation}
Let us introduce the scalars $\beta_i$, for $i=0,\cdots,q$ defined by $\beta_i=\frac{\gamma_i}{\gamma_q}$. In this case, we
have
\begin{equation}\label{betai}
\gamma_i= \frac{\beta_i}{\sum_{i=0}^{q} \beta_i} \; \text{for}\; i=0,\ldots,q-1 \; \text{and}\; \beta_q=1.
\end{equation} 
With these new variables, the system (\ref{eqgamma}) becomes
\begin{equation*}
\left\{
\begin{array}{ccccccccc}
\beta_0 (y_0,\Delta s_0) &+& \beta_1 (y_0,\Delta s_1) &+& \cdots &+& \beta_{q-1} (y_0,\Delta s_{q-1}) &=& - (y_0,\Delta s_{q})\\[6pt]
\vdots & & \vdots & &  & & \vdots & & \vdots \\[6pt]
\beta_0 (y_{q-1},\Delta s_0) &+& \beta_1 (y_{q-1},\Delta s_1) &+& \cdots &+& \beta_{q-1} (y_{q-1},\Delta s_{q-1}) &=& - (y_{q-1},\Delta s_{q}).
\end{array}
\right.
\end{equation*}
This system can be written as
\begin{equation}\label{linsyst}
   (Y_q^{T} \Delta S_q)\beta= -Y_q^{T} \Delta s_q, 
\end{equation}
where $\beta=(\beta_0,\cdots,\beta_{q-1})^{T}$ and $\Delta S_q=(\Delta s_0,\cdots, \Delta s_{q-1}).$\\
Now let us assume that the coefficients $\gamma_0,...,\gamma_q$ have been calculated and introduce the variables
\begin{equation}
    \begin{cases}
     \alpha_0=1-\gamma_0 \\
\alpha_{j}=\alpha_{j-1}-\gamma_{j} \,\text{for} \,j=1,\ldots,q-1.   
    \end{cases}
\end{equation}
In this case, the vector $t_q$ can be expressed as
\begin{equation}
    t_q= s_0 + \sum_{j=0}^{q-1} \alpha_j \Delta s_j = s_0 + \Delta S_q \alpha,
\end{equation}
where $\alpha=(\alpha_0,\cdots,\alpha_{q-1})^{T}$.
\begin{remark}
Note that to determine the unknowns $\gamma_i$, we must first calculate
the $\beta_i$ by solving the system of linear equations (\ref{linsyst}).
\end{remark}
\noindent We now define each extrapolation method and present the corresponding algorithms.\\
The RRE method is obtained by taking


\begin{equation}
    \eta_{i,j}=(\Delta^{2}s_{k+i}   , \Delta s_{k+j} ).
\end{equation}
$Y_q$ is therefore defined by
$$ Y_q= \Delta ^{2}S_{k,q}, $$
and the extrapolated approximation $t_{k,q}$ for the RRE method is given by
\begin{equation}
    t_{k,q}= s_{k}-\Delta S_{k,q}\Delta^{2} S_{k,q}^{+} \Delta s_{k},
    \label{eq:rreapproxim}
\end{equation}
where $\Delta^{2} S_{k,q}^{+}$ is the Moore-Penrose generalized inverse of $\Delta^{2} S_{k,q}$ defined by
$$  \Delta^{2} S_{k,q}^{+}=\left({\Delta^{2} S_{k,q}}^{T} \Delta^{2} S_{k,q} \right)^{-1} {\Delta^{2} S_{k,q}}^{T}. $$
In particular, the generalized residual of $t_{k,q}$ for the RRE method is defined as
\begin{equation}
    \tilde{r}(t_{k,q})= \Delta s_{k}-\Delta ^{2} S_{k,q}\Delta^{2} S_{k,q}^{+} \Delta s_{k},
    \label{eq:rreapproxim}
\end{equation}

\noindent An effective and reliable implementation for the reduced rank and the minimal polynomial extrapolation methods was provided by \cite{sidi1991efficient}. See \cite{brezinski1975generalisations,ford1988recursive} to examine other recursive methods suggested in the literature for implementing these polynomial vector extrapolation strategies. For more
details about the algorithms and convergence analysis of these methods, we refer to \cite{sidi1986convergence,jbilou2000vector,sidi1991efficient,jbilou1991some,jbilou1995analysis,sidi2017vector}.\\
Algorithm \ref{algo3} describes the RRE method \cite{sidi1991efficient,duminil2011reduced,duminil2015nonlinear,duminil2014fast}.

\begin{algorithm}[H]
\caption{The RRE method}\label{algo3}
\begin{algorithmic}[1]
\Require $\text{Vectors} \,s_{0}, s_{1},...,s_{q+1}$
\Ensure $t_{q}^{\text{RRE}}: \text{the RRE extrapolated  approximation} $
\State $ \text{Compute}\, \Delta s_{i}= s_{i+1}-s_{i},\,\text{for} \, i=0,1,...,q$
\State $\text{Set } \Delta U_{q+1}=[\Delta s_{0},\Delta s_{1},...,\Delta s_{q}]$
\State $\text{Compute the QR factorization of}\, \Delta U_{q+1}, \text{namely}, \,\Delta U_{q+1}=Q_{q+1}R_{q+1}$
\State $ \text{Solve the linear system}\, R_{q+1}^{T}R_{q+1}d=e, \,\text{where} \,d=[d_0,...,d_q]^{T} \,\text{and}\, e=[1,...,1]^{T}$
\State $ \text{Set}\, \lambda=(\sum_{i=0}^{q}d_{i})\, \text{and } \,\gamma_{j}= \frac{1}{\lambda}d_{i}, \,\text{for} \,i=0,...,q$
\State $ \text{Compute } \alpha=[\alpha_0,\alpha_1,...,\alpha_{q-1}]^{T}\, \text{where}\, \alpha_0=1-\gamma_0 \,\text{and } \,\alpha_{j}=\alpha_{j-1}-\gamma_{j} \,\text{for} \,j=1,...,q-1$
\State $ \text{Compute}\, t_{q}^{\text{RRE}}=s_0+ Q_{q}(R_{q}\alpha)$
\end{algorithmic}
\end{algorithm}

\noindent The scalars $\eta_{i,j}$ for the MPE method are defined by
\begin{equation}
    \eta_{i,j}=(\Delta s_{k+i}   , \Delta s_{k+j} ).
\end{equation}
$Y_q$ is therefore defined by
$$ Y_q= \Delta S_{k,q}, $$
and the extrapolated approximation $t_{k,q}$ is given by
\begin{equation}
    t_{k,q}= s_{k}-\Delta S_{k,q}\left({\Delta S_{k,q}}^{T} \Delta^{2} S_{k,q} \right)^{-1} {\Delta S_{k,q}}^{T} \Delta s_{k},
    \label{eq:rreapproxim}
\end{equation}

\noindent The MPE method \cite{sidi1991efficient} is implemented in full in Algorithm~\ref{algo4}.
\begin{algorithm}[H]
\caption{The MPE method}\label{algo4}
\begin{algorithmic}[1]
\Require $\text{Vectors} \,s_{0}, s_{1},...,s_{q+1}$
\Ensure $t_{q}^{\text{MPE}}: \text{the MPE extrapolated  approximation} $
\State $ \text{Compute}\, \Delta s_{i}= s_{i+1}-s_{i},\,\text{for} \, i=0,1,...,q$
\State $\text{Set } \Delta U_{q+1}=[\Delta s_{0},\Delta s_{1},...,\Delta s_{q}]$
\State Compute the QR factorization of $\Delta U_{q+1}$, namely, \,$\Delta U_{q+1}=Q_{q+1}R_{q+1} (\Delta U_{q}=Q_{q}R_{q}$ is contained in \, $\Delta U_{q+1}=Q_{q+1}R_{q+1})$
\State Solve the upper triangular linear system $R_{q}d=-r_q,$ \,\text{where} \,$d=[d_0,...,d_q]^{T}$ \,\text{and}\,$ r_{q}=[r_{0q},r_{1q},...,r_{q-1,q}]^{T}$
\State $ \text{Set}\,d_q=1 \text{and compute}\, \lambda=(\sum_{i=0}^{q}d_{i})\, \text{and set} \,\gamma_{j}= \frac{1}{\lambda}d_{i}, \,\text{for} \,i=0,...,q$
\State $ \text{Compute } \alpha=[\alpha_0,\alpha_1,...,\alpha_{q-1}]^{T}\, \text{where}\, \alpha_0=1-\gamma_0 \,\text{and } \,\alpha_{j}=\alpha_{j-1}-\gamma_{j} \,\text{for} \,j=1,...,q-1$
\State $ \text{Compute}\, t_{q}^{\text{MPE}}=s_0+ Q_{q}(R_{q}\alpha)$
\end{algorithmic}
\end{algorithm}

\subsection{The vector $\varepsilon$-algorithm}
Wynn \cite{wynn1956device} identified the $\epsilon$-algorithm as an effective way to apply the Shanks transformation \cite{brezinski1975generalisations}. It is a method developed to speed up the convergence of a sequence $S=(s_0,s_1,s_2,\ldots,s_j\in\mathbb{R})$ generated by a linear or nonlinear fixed-point iterative method.\\
We give a description on the implementation of the vector $\epsilon$-algorithm of Wynn \cite{wynn1956device}.\\
\textit{Initialization:} \\
For $j = 0,1,2,\ldots$
\begin{align}
    \varepsilon^{(j+1)}_{-1} &:= 0,  \\[6pt]
    \varepsilon^{(j)}_{0} &:= s_{j}. 
\end{align}

\textit{Iteration:}\\
\quad For $k,j = 0,1,2,\ldots$

\begin{equation}\label{vectorepsilonapprox}
    \varepsilon^{(j)}_{k+1} 
    := \varepsilon^{(j+1)}_{k-1} 
      + \big[\, \varepsilon^{(j+1)}_{k} - \varepsilon^{(j)}_{k} \,\big]^{-1},
\end{equation}
where the inverse of a vector $\mathrm{v}=(v_1, v_2,\ldots, v_n)^{T}\in \mathbb{R}^{N}$ is given by
\begin{equation}
	\mathrm{v}^{-1}:= \frac{\mathrm{v}}{\mathrm{v} \cdot \mathrm{v}}.
\end{equation}
For given values of $J , K$, the extrapolated approximation is given by $\epsilon_{2K}^{(J)}$.\\
With this definition (\ref{vectorepsilonapprox}) reduces to the scalar algorithm for $n=1$. \\

\begin{figure}[h!]
\[
\begin{array}{ccccccccc}
 & \varepsilon^{(0)}_{0} &    &   &       &  &       &  \\[8pt]
\varepsilon^{(1)}_{-1} &    & \varepsilon^{(0)}_{1} &        &  &        &  &  \\[8pt]
& \varepsilon^{(1)}_{0} &        & \varepsilon^{(0)}_{2} &      &   &  &      \\[8pt]
\varepsilon^{(2)}_{-1} &      & \varepsilon^{(1)}_{1} & &  \cdot    &       &        &      \\[8pt]
 &  \varepsilon^{(2)}_{0} &       & \varepsilon^{(1)}_{2}  &       &  \cdot &        & \\[8pt]       
\varepsilon^{(3)}_{-1}     &   &    \varepsilon^{(2)}_{1}    &  &    \cdot   & &     \cdot   &  \\[8pt] 
   &  \cdot &      & \cdot &       & \cdot  &      & \cdot \\[8pt] 
\end{array}
\]
\caption{The vector $\varepsilon$-table.}
\end{figure}

\section{Solution of nonlinear least-squares problems} \label{sec3}
\subsection{Iterative methods}
The Gauss--Newton method \cite{buccini2023efficient} is one among many iterative techniques \cite{dennis1996numerical,ruhe1979accelerated,pes2020minimal,pes2022doubly} developed for solving nonlinear least-squares problems as well as other similar problems. In \cite{buccini2023efficient}, the authors proposed a method based on projection onto generalized Krylov subspaces. More specifically, they presented an efficient implementation of the Gauss--Newton method through the use of Generalized Krylov Subspaces (GKS), resulting in the GNGKS method (as referred to as GNKS in this work) for solving large nonlinear least-squares problems.\\
Consider the nonlinear least-squares problem 
\begin{equation}\label{Lstsq}
  \arg \min _{\mathrm{x}\in \mathbb{R}^n} g(\mathrm{x})= \arg \min _{\mathrm{x}\in \mathbb{R}^n} {\Vert \mathrm{y} -f(\mathrm{x}) \Vert} ^{2},
 \end{equation}
with a nonlinear differentiable function $f$ and available data $\mathrm{y}\in \mathbb{R}^n$.\\
For an initial guess $\mathrm{x}_0\in \mathbb{R}^n$, the gradient descent iterations for solving the least-squares problem~(\ref{Lstsq}) are given by
\begin{equation}\label{gradientdescent}
\mathrm{x}_{k+1} = \mathrm{x}_k - \tau_k \nabla g(\mathrm{x}_k), \;\;\;\;\;\;\;\;\;\;\;\;\;\;\;k=0,1,2,\ldots
\end{equation}
where $\nabla g$ denotes the gradient of the nonlinear function $g$, and $\tau_k$ is determined by a line search rule such as the Armijo condition (\ref{Armijo-Goldstein}) with backtracking (see \cite{buccini2023efficient})
\begin{equation}\label{Armijo-Goldstein}
g(\mathrm{x}_{k+1}) \leq g(\mathrm{x}_{k}) - \omega \tau_{k} \langle \nabla g(\mathrm{x}_k), \nabla g(\mathrm{x}_k) \rangle,
\end{equation}
with constant $\omega$. 

A general preconditioned form of the gradient descent iterative method is given by 
\begin{equation}\label{preconditgradientdescent}
\mathrm{x}_{k+1} = \mathrm{x}_k - \tau_k H_{k}^{-1}\nabla g(\mathrm{x}_k), \;\;\;\;\;\;\;\;\;\;\;\;\;\;\;k=0,1,2,\ldots
\end{equation}
where the constant $\tau_k$ satisfies
\begin{equation}\label{preconditArmijo-Goldstein}
g(\mathrm{x}_{k+1}) \leq g(\mathrm{x}_{k}) - \omega \tau_{k} \langle  H_{k}^{-1} \nabla g(\mathrm{x}_k), \nabla g(\mathrm{x}_k) \rangle.
\end{equation}
We use a line search to find an appropriate $\tau_k$. Let $\tau^{(0)}_k > 0$ be sufficiently large. We set $\tau_k= \tau^{(0)}_k$ if $\tau^{(0)}_k$ meets condition (\ref{preconditArmijo-Goldstein}); if not, we let
$$\tau^{(1)}_k=\frac{\tau^{(0)}_k}{2}.$$
This procedure is repeated until an appropriate step size $\tau_k^{(j)}$ satisfying~(\ref{preconditArmijo-Goldstein}) is found. We then set $\tau_k = \tau_k^{(j)}$.

If the matrix $H_k$ is the identity, we obtain the classical gradient descent iterations (\ref{gradientdescent}).\\
We choose a preconditioner matrix $H_k = \mathrm{diag}(J_f)$, which is the diagonal of the Jacobian matrix of $f$, and we denote by PGD the iterative process generated by the preconditioned gradient descent (Algorithm \ref{algooo1}). Moreover, if we choose $H_k = \mathrm{diag}(J_f^{T} J_f)$, where $J_f^{T}$ is the transpose of the Jacobian matrix of $f$, we obtain the scaled gradient descent (Algorithm \ref{algooo2}), which is denoted by SGD. It should be noted that the importance of this SGD method is particularly evident when the Jacobian matrix $J_f$ is not square.


\begin{algorithm}[H]
\caption{Preconditioned Gradient Descent (PGD) algorithm}\label{algooo1}
\begin{algorithmic}[1]
\Require Initial guess $\mathrm{x}_0$, $\text{Itermax}$, a constant $\omega_1$
\Ensure  $\mathrm{x}_k$
\State $k \gets 0$
\While{$k \leq \text{Itermax}$ \text{and not converged} }
    \State Determine the preconditioner $H_k = \mathrm{diag}(J_f(\mathrm{x}_k))$
    \State Calculate the direction of the descent $d_k = - H_k^{-1} \nabla g(\mathrm{x}_k)$
    \State Choose $\tau_k$ satisfying $g(\mathrm{x}_{k+1}) \leq g(\mathrm{x}_{k}) - \omega_{1}\,\tau_{k} \langle H_{k}^{-1}\nabla g(\mathrm{x}_k), \nabla g(\mathrm{x}_k) \rangle $
    \State Update the iterate $\mathrm{x}_{k+1} = \mathrm{x}_k + \tau_k d_k$
    \State $k \gets k+1$
\EndWhile
\end{algorithmic}
\end{algorithm}

\begin{algorithm}[H]
\caption{Scaled Gradient Descent (SGD) algorithm}\label{algooo2}
\begin{algorithmic}[1]
\Require Initial guess $\mathrm{x}_0$, $\text{Itermax}$, a constant $\omega_2$
\Ensure  $\mathrm{x}_k$
\State $k \gets 0$
\While{$k \leq \text{Itermax}$ \text{and not converged} }
    \State Determine $H_k = \mathrm{diag}(J_f^{T}(\mathrm{x}_k)J_f(\mathrm{x}_k))$
    \State Calculate   $d_k = - H_k^{-1} \nabla g(\mathrm{x}_k)$
    \State Choose $\tau_k$ by $g(\mathrm{x}_{k+1}) \leq g(\mathrm{x}_{k}) - \omega_{2}\,\tau_{k} \langle H_{k}^{-1}\nabla g(\mathrm{x}_k), \nabla g(\mathrm{x}_k) \rangle $
    \State Update the iterate $\mathrm{x}_{k+1} = \mathrm{x}_k + \tau_k d_k$     
    \State $k \gets k+1$
\EndWhile
\end{algorithmic}
\end{algorithm}

\subsection{An extrapolation-based method for nonlinear least-squares problems }
The methods given in Algorithms~\ref{algo3} and~\ref{algo4} become increasingly expensive as \( q \) grows. For the complete forms of the RRE and MPE methods, the computational effort increases quadratically with the number of iterations \( q \), while the storage cost increases linearly. The computational cost, in terms of arithmetic operations (additions and multiplications), and the memory requirements for computing the approximation \( t_{k,q} \) using RRE and MPE are, respectively, \( 2Nq^2 \) operations and \( q+1 \) stored vectors, where \( N \) denotes the problem dimension. Moreover, when applied to the solution of linear systems, both RRE and MPE require \( q+1 \) matrix--vector products per extrapolation cycle. A practical way to keep both storage and computational costs at a reasonable level is to restart these algorithms periodically~\cite{mouhssine2025vector,mouhssine2025enhancing}. We denote the restarted (cyclic) versions of the extrapolation techniques by RRE(\( q \)) and MPE(\( q \)). Similarly, the restarted version of the vector \( \varepsilon \)-algorithm is denoted by VEA(\( q \)).

The sequence $(\mathrm x_k)_k$ (\ref{preconditgradientdescent}) generated by the preconditioned gradient descent method is then used to construct the vector sequence $(s_k)_k$ (\ref{fixedpoint}), with $s_k = \mathrm x_k$, to which vector extrapolation techniques are applied in order to accelerate convergence.

Algorithm \ref{algoOOOO5} describes the restarted approach that combines the polynomial extrapolation methods like RRE, and MPE with the preconditioned gradient descent iterations (PGD and SGD).

\begin{algorithm}[H]
\caption{A polynomial extrapolation algorithm for preconditioned gradient descent}\label{algoOOOO5}
\begin{algorithmic}[1]
\Require $\mathrm x_0$, $\text{fixed restart number}\,q$, $\text{Itermax}$
\Ensure $t_{q}: \text{the extrapolated  approximation} $
\State $k \gets 0$
\While{$k \leq \text{Itermax and not converged}$}
\While{$i \leq q$}
 \State $ \text{Compute }\, \mathrm x_{i+1}\,\text{by Algorithm} \,\ref{algooo1} \, \text{or}\, \text{Algorithm}\,\ref{algooo2}$
 \State $i \gets i + 1$
\EndWhile
\State  $\text{Calculate}\,t_{q}\, \text{using Algorithm}\, \ref{algo3}\,  \text{or}\, \text{Algorithm}\, \ref{algo4} $
\If{$t_q\, \text{satisfies accuracy test}$}
        \State $\textbf{return}\, t_q $
\Else
        \State $\mathrm x_0 \gets t_q$
        \State $k \gets k+ 1$
\EndIf
\EndWhile

\end{algorithmic}
\end{algorithm}

\noindent Moreover, the restarted version of the VEA is described in Algorithm \ref{algogchhcjcs} below and is applied to the iterations generated by the preconditioned gradient descent method.
 
\begin{algorithm}[H]
\caption{Vector Epsilon Algorithm (VEA) for preconditioned gradient descent}\label{algogchhcjcs}
\begin{algorithmic}[1]
\Require $\text{Given an initial guess} \,\mathrm x_0$ \text{and} $q$: the restart number, $\text{Itermax}$
\Ensure $\varepsilon_{2q}: \text{the VEA extrapolated  approximation} $
\State $k \gets 0$
\While{$k \leq \text{Itermax and not converged}$}
\While{$i \leq 2q -1$}
 \State $ \text{Compute }\, \mathrm x_{i+1}  \text{by Algorithm} \,\ref{algooo1} \, \text{or}\, \text{Algorithm}\,\ref{algooo2}$
 \State $i \gets i + 1$
\EndWhile
\State  $\text{Compute}\,\varepsilon_{2q}\, \text{using the VEA process}$
\If{$\varepsilon_{2q}\, \text{satisfies accuracy test}$}
        \State $\textbf{return}\, \varepsilon_{2q} $
\Else
        \State $\mathrm x_0 \gets \varepsilon_{2q}$
        \State $k \gets k+ 1$
\EndIf
\EndWhile

\end{algorithmic}
\end{algorithm}

\begin{remark}
It should be highlighted that in contrast to polynomial extrapolation techniques (RRE and MPE), which require the generation of $q+2$ vectors (in particular, the preconditioned gradient descent iterations), including the initial vector $\mathrm x_0$, in the vector epsilon algorithm (VEA), $2q+1$ vectors, including the initial guess, must be generated in each cycle.   
\end{remark}

\section{Numerical experiments}   \label{sec4}
To solve nonlinear problems, we numerically investigate the effectiveness of vector extrapolation techniques combined with gradient descent methods and their preconditioned variants. The proposed approaches are compared with the Gauss--Newton method based on generalized Krylov subspaces \cite{buccini2023efficient}.

\subsection{Numerical study of the nonlinear Bratu problem}
\subsubsection{Model problem and numerical setting}
Consider the model problem (\ref{eq2dbratu})
\begin{equation}
 \begin{cases}
-\text{div}(\nabla x(s,t)) + \alpha \frac{\partial x_ (s,t)}{\partial s}+ \lambda e^{x(s,t)}& = y \;\;\;\;\;\;\;  \text{in} \;\;  \Omega, \\
\;\;\;\;\;  x(s,t) &= 0 \;\;\;\;\;\;\; \text{on} \;\; \Gamma=\partial \Omega,
\end{cases}
\label{eq2dbratu}
\end{equation}
where $\Omega=[-3, 3] \times [-3, 3]$ represents the domain of our problem. In both the $s$ and $t$ directions, we discretize this equation (\ref{eq2dbratu}) using uniform grids with different grid spacings with $n$ discretization points. The Laplacian is approximated using a typical second-order finite difference operator. It is easy to see that the discretized problem can be rewritten as a minimization problem, namely  
\begin{equation*}\label{lstsq}
  \arg \min _{\mathrm x\in \mathbb{R}^{n^{2}}} g(\mathrm x)= \arg \min _{\mathrm x\in \mathbb{R}^{n^{2}}} {\Vert \mathrm{y} -f(\mathrm x) \Vert} ^{2},
 \end{equation*}
where 
\[
f(\mathrm{x}) = L \mathrm{x} + \alpha D \mathrm{x} + \lambda e^{\mathrm{x}},
\]
and $\mathrm{x}$ denotes the discretized solution vector.  
The matrices $L$ and $D$ represent discrete differential operators. Specifically, we set
\[
L = L_1 \otimes I + I \otimes L_1, \qquad 
D = D_1 \otimes I,
\]
where $\otimes$ denotes the Kronecker product, and
\[
L_1 =
\begin{bmatrix}
	2 & -1 &  &   \\
	-1 & 2 & -1 &  \\
	& \ddots & \ddots & \ddots \\
	&  & -1 & 2
\end{bmatrix}
\in \mathbb{R}^{n\times n}, 
\quad
D_1 =
\begin{bmatrix}
	-1 & 1 &  &   \\
	& \ddots & \ddots &  \\
	&  & -1 & 1
\end{bmatrix}
\in \mathbb{R}^{n\times n}.
\]
The operator $L$ corresponds to the discrete Laplacian, 
while $D$ represents the first-derivative operator in the $s$-direction.  
Once the discretized function $\mathrm{x}_{\text{true}}$ is sampled (by computing the sampling of the function $x$), 
we calculate $\mathrm{y}$ by
\[
\mathrm{y} = L \mathrm{x}_{\text{true}} + \alpha D \mathrm{x}_{\text{true}} + 
\lambda e^{\mathrm{x}_{\text{true}}}.
\]
In our numerical tests, we consider the exact solution 
$x(s,t) = e^{-10(s^{2}+t^{2})}$, 
sampled on a uniform $100\times100$ grid over the square domain $[-3,3]^2$. 
Therefore, in (\ref{Lstsq}), $f$ is defined as 
$f:\mathbb{R}^{10^{4}} \rightarrow \mathbb{R}^{10^{4}}$. 
The Jacobian of $f$ at a point $\mathrm{x}$ is then given by
\[
J_{f} = L + \alpha \cdot D + \lambda\,\cdot \mathrm{diag}(e^{\mathrm{x}}).
\]

\subsubsection{Numerical study of the extended Bratu problem}
We illustrate the application of polynomial extrapolation methods to gradient descent iterations, with and without preconditioning, and assess the performance of the resulting hybrid approaches for the solution of the two-dimensional nonlinear Bratu problem~(\ref{eq2dbratu}).

\begin{figure}[htbp]
    \centering
     \includegraphics[width=4.4cm,height=5cm]{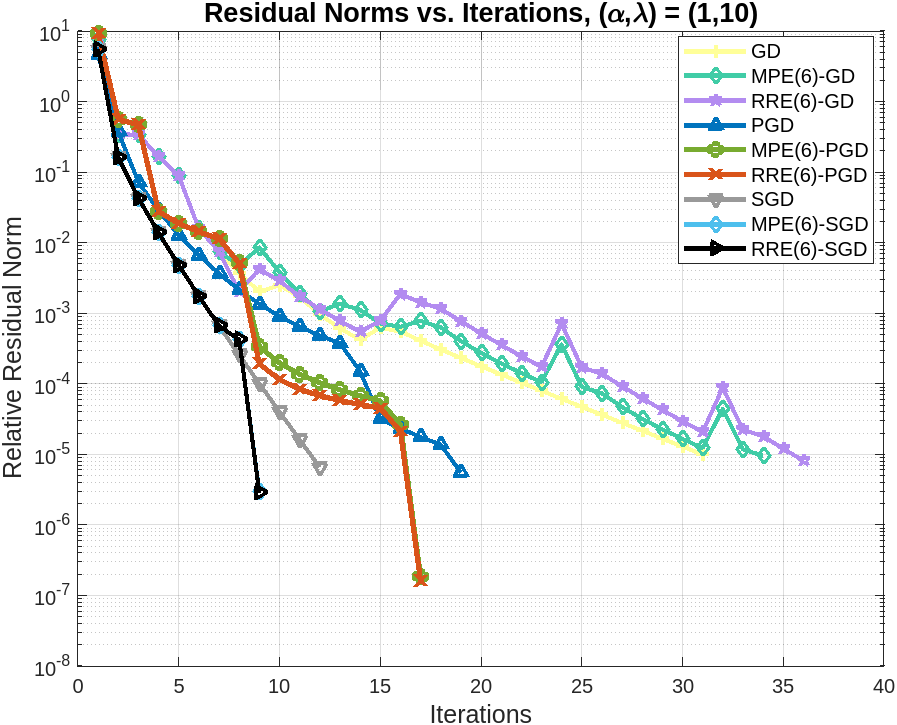}
   \includegraphics[width=4.4cm,height=5cm]{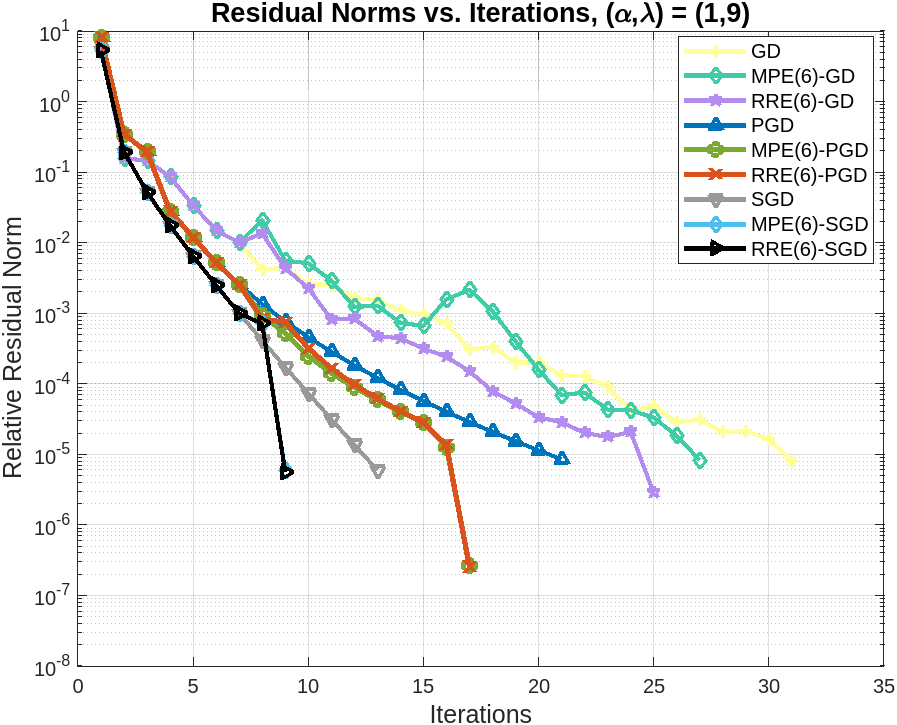}
   \includegraphics[width=4.4cm,height=5cm]{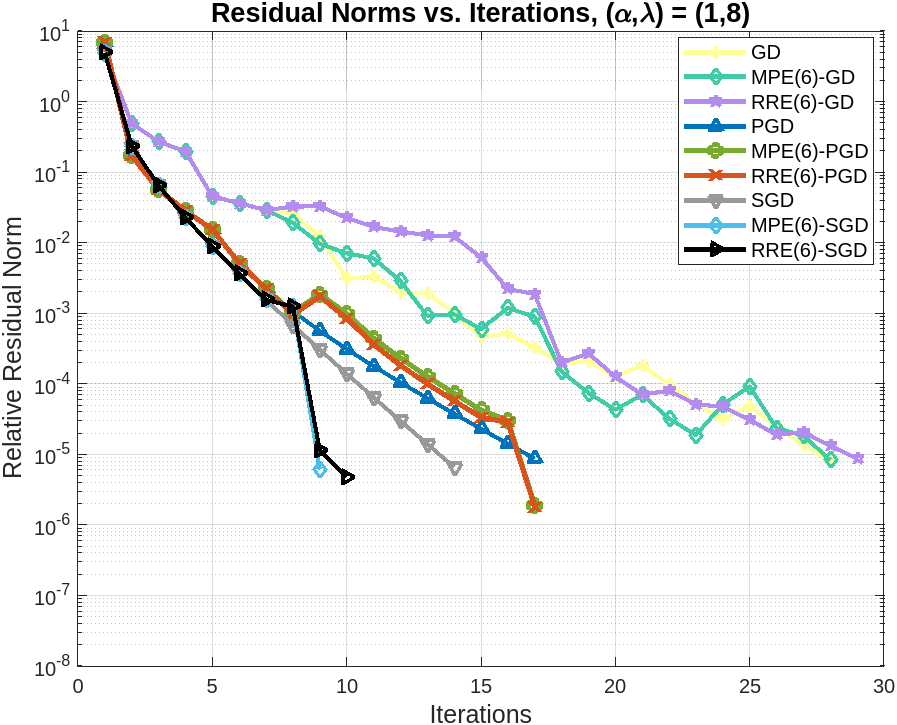}
  
   \caption{Performance of gradient descent and its preconditioned versions with and without polynomial extrapolation.}
   \label{figbenchmark_residual}
\end{figure}

First, Figure~\ref{figbenchmark_residual} illustrates the performance of the proposed hybrid solvers for selected values of the parameters $\alpha$ and $\lambda$ in a well-conditioned regime of the Bratu equation, corresponding to small values of $\alpha$ and relatively large values of $\lambda$. For all tests, a $10^{2}\times 10^{2}$ interior grid is used, and the iterations are stopped when the relative residual norm, defined as the $\ell_2$-norm of the difference between two successive iterates,$
\frac{\Vert \mathrm{x}_{k+1}-\mathrm{x}_{k} \Vert_2}{\Vert \mathrm{x}_{k} \Vert_2},$
falls below $10^{-5}$.

The baseline gradient-based methods without acceleration are denoted by GD, PGD, and SGD. In all numerical experiments, the parameters $\omega_1=10^{-4}$ and $\omega_2=0.5$ are used for the PGD and SGD methods, respectively. In addition, we assess the efficiency of polynomial extrapolation techniques, namely RRE and MPE, for accelerating the convergence of gradient descent and preconditioned gradient descent iterations. The corresponding accelerated methods are denoted by RRE($q$)-GD, MPE($q$)-GD, RRE($q$)-PGD, MPE($q$)-PGD, RRE($q$)-SGD, and MPE($q$)-SGD.

\noindent We observe that neither RRE nor MPE significantly improves the convergence of the classical gradient descent (GD) iterations. Even when varying the number of elements $q$ used to construct the extrapolated approximation, the resulting convergence remains slow and unstable when extrapolation is applied directly to GD. For this reason, we consider preconditioned variants of gradient descent, namely PGD and SGD, which already provide a substantial acceleration compared to the standard GD method for solving~(\ref{eq2dbratu}). When combined with these preconditioned schemes, polynomial extrapolation techniques yield an additional acceleration, as illustrated in Figure~\ref{figbenchmark_residual}.

This behavior highlights the effectiveness of extrapolation techniques when applied to preconditioned gradient-based methods rather than to their unpreconditioned versions.

\begin{figure}[htbp]
    \centering
     \includegraphics[width=4.4cm,height=5cm]{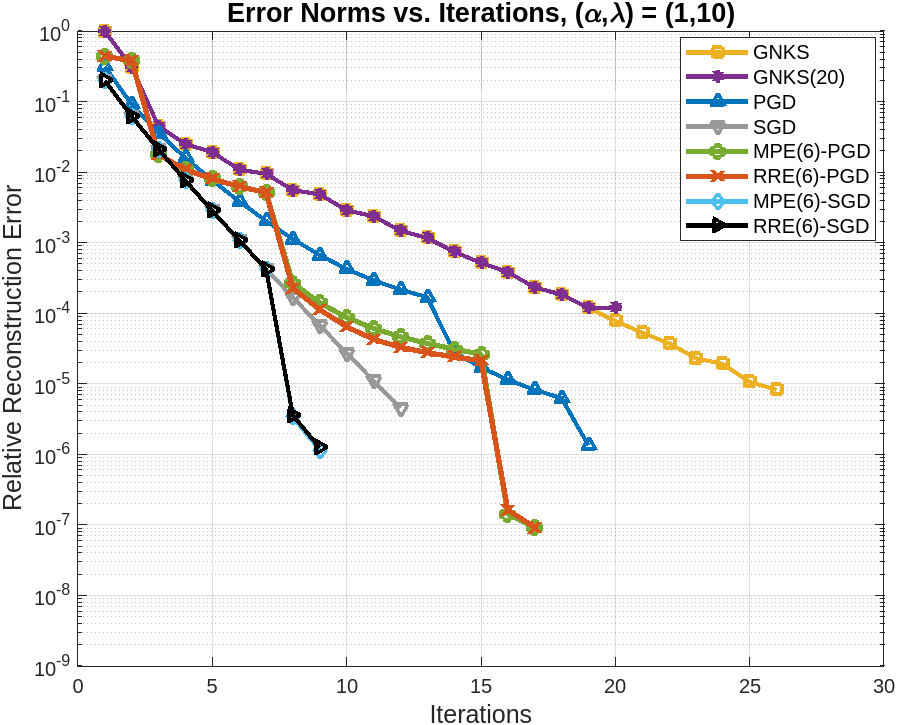}
   \includegraphics[width=4.4cm,height=5cm]{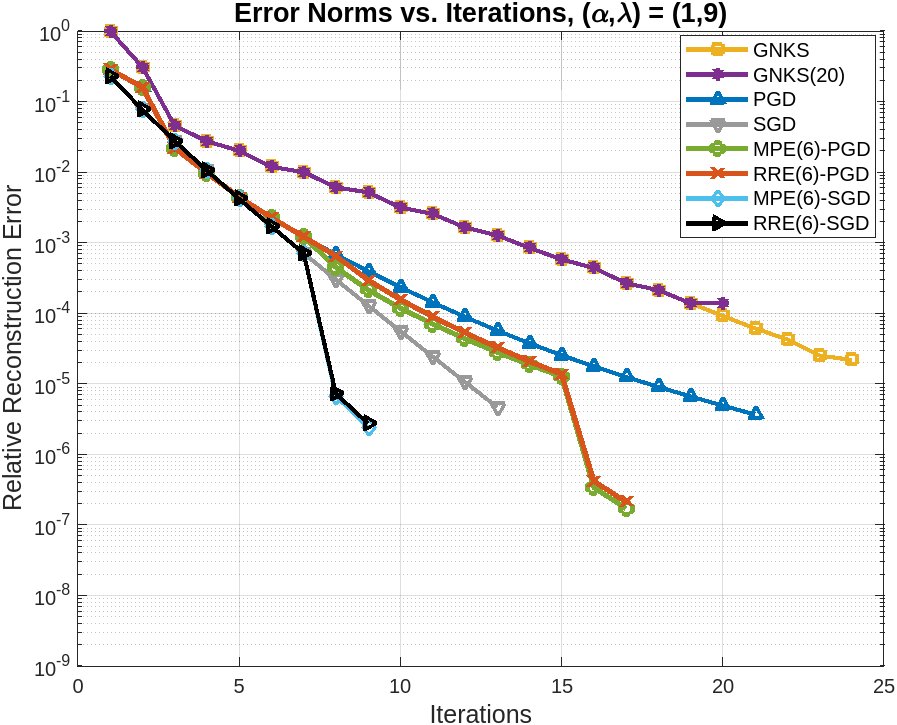}
   \includegraphics[width=4.4cm,height=5cm]{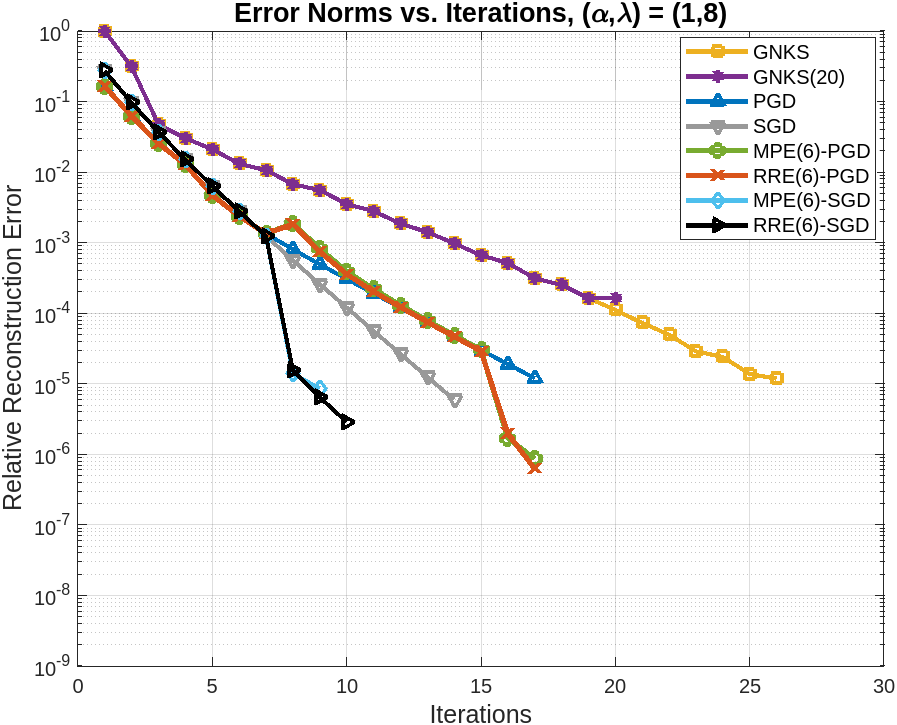}
  
   \caption{Convergence of gradient descent and its preconditioned versions with and without polynomial extrapolation in comparison with Gauss--Newton using generalized Krylov subspaces.}
   \label{figbenchmark}
\end{figure}

Figure~\ref{figbenchmark} shows the relative reconstruction error,$
\frac{\Vert \mathrm{x}-\mathrm{x}_{\text{true}} \Vert_{2}}{\Vert \mathrm{x}_{\text{true}} \Vert_{2}},$
as a function of the iteration number, where $\mathrm{x}_{\text{true}}$ denotes the reference solution. For comparison purposes, we also report the results obtained with the Gauss--Newton method based on generalized Krylov subspaces (GNKS) \cite{buccini2023efficient} and its restarted variant GNKS($20$), with restart parameter $q=20$.

\noindent We observe that the restarted polynomial extrapolation methods with $q=6$, namely MPE(6)-PGD, RRE(6)-PGD, MPE(6)-SGD, and RRE(6)-SGD, achieve a lower relative $L_2$ error in fewer iterations than the PGD, SGD, GNKS, and GNKS(20) methods. Moreover, Figure~\ref{figbenchmark} indicates that the MPE-PGD and RRE-PGD approaches yield more accurate and efficient approximations of~(\ref{eq2dbratu}) than their SGD-based versions, namely MPE-SGD and RRE-SGD.

\begin{figure}[h!]
    \centering

    \begin{subfigure}[t]{0.48\textwidth}
        \centering
        \includegraphics[width=6cm,height=6cm]{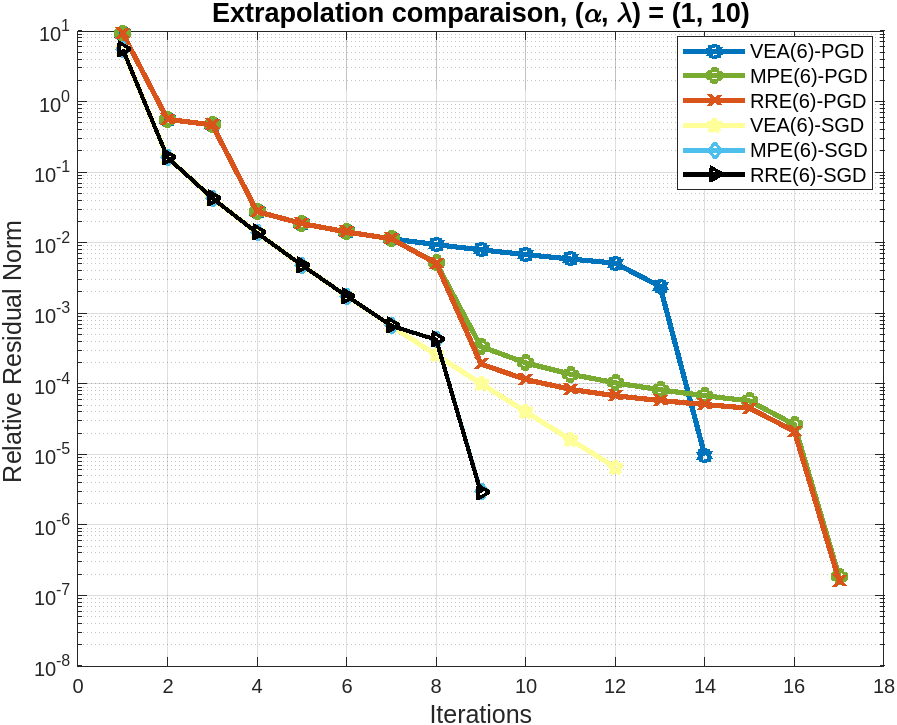}
        \caption{Residual norms}
        \label{fig:extrap_residual}
    \end{subfigure}
    \hfill
    \begin{subfigure}[t]{0.48\textwidth}
        \centering
        \includegraphics[width=6cm,height=6cm]{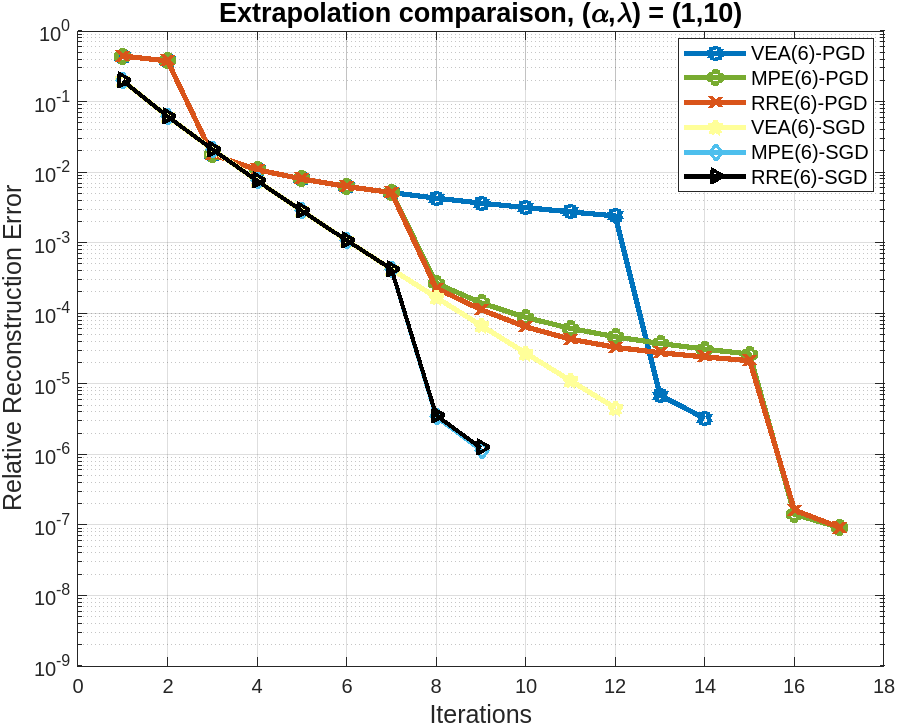}
        \caption{Error norms}
        \label{fig:extrap_error}
    \end{subfigure}

    \caption{Efficiency comparison between MPE, RRE, and VEA for the acceleration of the preconditioned gradient descent.}
    \label{fig:extrap_comparison}
\end{figure}

Figure~\ref{fig:extrap_comparison} illustrates the convergence behavior of vector extrapolation methods, including polynomial-type techniques such as RRE and MPE, as well as the vector $\epsilon$-algorithm (VEA), when applied to accelerate the iterations of preconditioned gradient methods, namely PGD and SGD. On the one hand, these methods are effective in improving the convergence rate. On the other hand, since we focus on restarted extrapolation strategies, polynomial extrapolation techniques are preferred, as they require storing fewer vectors. More precisely, they only require $q+2$ vectors, including the initial iterate, from the PGD or SGD iterations at each cycle to construct the extrapolated approximation, whereas the VEA method requires storing $2q+1$ vectors, including the initial guess.

\noindent Table~\ref{table_bratu} reports a benchmarking study for different values of the parameters $\alpha$ and $\lambda$ in the well-conditioned regime of~(\ref{eq2dbratu}), comparing the GNKS, GNKS(20), MPE-PGD, RRE-PGD, MPE-SGD, and RRE-SGD methods. The table provides the number of iterations required to reach the prescribed tolerance, together with the relative reconstruction error and the CPU time in seconds.

The results in Table~\ref{table_bratu} show that the MPE-PGD, RRE-PGD, MPE-SGD, and RRE-SGD methods achieve better convergence in terms of relative reconstruction error, producing more accurate approximate solutions than GNKS and GNKS(20) in fewer iterations and with comparable CPU times. In several cases, the MPE-PGD and RRE-PGD techniques further outperform their SGD-based methods in terms of both relative error and CPU time. Moreover, the convergence behavior of MPE-PGD and RRE-PGD is very similar. Consequently, we focus on the MPE-PGD and RRE-PGD techniques for solving the Bratu problem~(\ref{eq2dbratu}).

\begin{table}[htbp]
	 \caption{The Bratu equation for fixed grid size with $n=10^{2}$ on each direction: Iteration count, relative reconstruction error (RE), and CPU time in seconds for various values of $\alpha$ and $\lambda $ using several approaches.}
  \label{table_bratu}	
  \centering
     \begin{tabular}{|c|c|c|c|c|}
        \hline
      $\lambda$  & Method  & Iter & Relative Error (RE) &  CPU(s)  \\
        \hline
         \multicolumn{5}{|c|}{$\alpha=1$} \\
        \hline
     &  GNKS & 26 & 8.20e-06  & 0.35 \\ 
   10 &  GNKS(20) & 20 & 1.19e-04  & 0.11  \\ 
  & RRE(6)-PGD & 17 & \textbf{9.26e-08}  & 0.10 \\
  & MPE(6)-PGD & 17 & \textbf{9.20e-08} & 0.10  \\
   & RRE(6)-SGD & 9  & 1.23e-06  & 0.08 \\ 
  & MPE(6)-SGD   & 9 & 1.13e-06 & 0.07  \\ 
       \hline 
     &  GNKS & 24 & 2.17e-05  & 0.11 \\ 
   9 &  GNKS(20) & 20 & 1.38e-04  & 0.07  \\ 
    & RRE(6)-PGD& 17 & \textbf{2.11e-07}  & 0.09 \\
  & MPE(6)-PGD   & 17 &  \textbf{1.70e-07} & 0.08  \\
   & RRE(6)-SGD & 9 & 2.72e-06 & 0.07 \\ 
  & MPE(6)-SGD  & 9 & 2.43e-06 & 0.07 \\  
\hline
 &  GNKS & 26 & 1.18e-05  & 0.45 \\ 
   8 &  GNKS(20) & 20 & 1.64e-04  & 0.10  \\ 
    & RRE(6)-PGD & 17 & \textbf{6.26e-07}  & 0.07 \\
  &  MPE(6)-PGD   & 17 &  \textbf{8.63e-07} & 0.07  \\
   & RRE(6)-SGD  & 10 & 2.87e-06 & 0.08 \\  
  & MPE(6)-SGD   & 9 & 8.61e-06 & 0.07 \\  
\hline
 &  GNKS & 26 & 1.34e-05  & 0.08 \\ 
  7 &  GNKS(20) & 20 & 1.99e-04  & 0.05  \\ 
    & RRE(6)-PGD & 15 & \textbf{4.67e-06}  & 0.08 \\
  & MPE(6)-PGD  & 15 &  \textbf{3.56e-06} & 0.08  \\
   & RRE(6)-SGD & 11 & 5.86e-06 & 0.09 \\  
  & MPE(6)-SGD  & 11 &  4.39e-06 &  0.08\\  
\hline
 &  GNKS & 26 & 2.12e-05  & 0.08 \\ 
   6 &  GNKS(20) & 20 & 2.52e-04  & 0.06  \\ 
    & RRE(5)-PGD & 16 & \textbf{8.14e-06}  & 0.05\\
  & MPE(5)-PGD   & 17 &  \textbf{1.41e-06} & 0.05  \\
   & RRE(6)-SGD & 11 & 2.00e-05 & 0.09 \\   
  & MPE(6)-SGD & 12 & 8.04e-06 & 0.10 \\  

\hline

     \end{tabular}

\end{table}

\begin{table}[htbp]
	 \caption*{\textbf{Table 1} continued.}
  \centering
\begin{tabular}{|c|c|c|c|c|}
        \hline
      $\lambda$  & Method  & Iter & Relative Error (RE) &  CPU(s)  \\
      \hline
   
              \multicolumn{5}{|c|}{$\alpha=2$} \\

                \hline
     &  GNKS & 24 & 1.73e-05  & 0.05 \\ 
   10 &  GNKS(20) & 20 & 8.99e-05  & 0.04  \\ 
    & RRE(6)-PGD & 17 & \textbf{6.26e-08}  & 0.06 \\
  & MPE(6)-PGD  & 17 &  \textbf{7.98e-08} & 0.05  \\
   & RRE(6)-SGD & 9 & 4.66e-07 & 0.07 \\ 
  & MPE(6)-SGD & 9 & 4.37e-07 & 0.08 \\  
       \hline 
                &  GNKS & 24 & 2.06e-05  & 0.05 \\ 
   9 &  GNKS(20) & 20 & 1.06e-04  & 0.03  \\ 
    & RRE(5)-PGD & 15 & \textbf{4.19e-07}  & 0.06 \\
  & MPE(5)-PGD   & 15 &  \textbf{3.46e-07} & 0.05  \\
   & RRE(6)-SGD & 9 & 8.54e-07 & 0.07  \\  
  & MPE(6)-SGD & 9 & 7.70e-07 & 0.07 \\ 
             \hline 
         &  GNKS & 24 & 2.40e-05  & 0.06 \\ 
  8 &  GNKS(20) & 20 & 1.26e-04  & 0.04  \\ 
    & RRE(7)-PGD & 19 & \textbf{7.53e-08}  & 0.06 \\
  & MPE(7)-PGD  & 19 &  \textbf{9.46e-08} & 0.05  \\
  & RRE(6)-SGD & 9 & 1.53e-06 & 0.08 \\  
  & MPE(6)-SGD & 9 & 1.33e-06 & 0.07 \\  
       \hline 
         &  GNKS & 24 & 2.63e-05  & 0.06 \\ 
  7 &  GNKS(20) & 20 & 1.52e-04  & 0.04  \\ 
    & RRE(10)-PGD & 17 & \textbf{5.73e-07}  & 0.05 \\
  & MPE(6)-PGD   & 18 &  \textbf{7.31e-07} & 0.05  \\
&  RRE(6)-SGD &9  & 3.59e-06 & 0.08 \\  
 &  MPE(6)-SGD & 9 & 3.12e-06 & 0.07 \\   
 
       \hline 
         &  GNKS & 26 & 1.63e-05  & 0.06 \\ 
   6 &  GNKS(20) & 20 & 1.85e-04  & 0.04  \\ 
    & RRE(8)-PGD & 15 & \textbf{2.59e-06 } & 0.06 \\
  & MPE(6)-PGD   & 18 &  \textbf{6.87e-07} & 0.06  \\
  & RRE(6)-SGD & 10 & 4.80e-06 & 0.08  \\  
  & MPE(6)-SGD & 10 & 4.02e-06 & 0.09 \\   
       \hline 
     \end{tabular}
\end{table}

\begin{table}[htbp]
	\caption*{\textbf{Table 1} continued.}
  \label{new_new_table_bratu}
  \centering
\begin{tabular}{|c|c|c|c|c|}
        \hline
      $\lambda$  & Method  & Iter & Relative Error (RE) &  CPU(s)  \\
      \hline     
       \multicolumn{5}{|c|}{$\alpha=3$} \\

      \hline
     &  GNKS & 24 & 2.37e-05  & 0.05 \\ 
  10 &  GNKS(20) & 20 & 1.24e-04  & 0.04  \\ 
 & RRE(6)-PGD & 17 & \textbf{2.51e-07} & 0.06 \\
& MPE(6)-PGD  & 17 & \textbf{2.12e-07} & 0.06  \\
  & RRE(6)-SGD & 8 & 2.48e-05 & 0.07 \\   
  & MPE(6)-SGD & 8 & 2.48e-05 & 0.07 \\ 
       \hline 
       &  GNKS & 22 & 5.47e-05  & 0.05 \\ 
  9 &  GNKS(20) & 20 & 1.29e-04  & 0.04  \\ 
    & RRE(7)-PGD & 18 & \textbf{3.36e-07}  & 0.06 \\
  & MPE(6)-PGD   & 17 &  \textbf{8.67e-07} & 0.05  \\
  & RRE(6)-SGD & 11 & 6.26e-06 & 0.10  \\  
  & MPE(6)-SGD & 11 & 6.23e-06 & 0.09 \\  
       \hline
       &  GNKS & 26 & 2.12e-05  & 0.06 \\ 
   8 &  GNKS(20) & 20 & 2.03e-04  & 0.04  \\
    & RRE(4)-PGD & 15 & \textbf{6.11e-06}  & 0.05 \\
  & MPE(3)-PGD   & 17 &  \textbf{3.90e-07} & 0.05  \\
   & RRE(6)-SGD & 12 & 7.38e-06 & 0.10   \\ 
  & MPE(6)-SGD & 12 & 7.35e-06 & 0.10 \\   
       \hline 
     \multicolumn{5}{|c|}{$\alpha=4$} \\

        \hline
  10   &  GNKS & 16 &  9.77e-04& 0.06  \\ 
    & RRE(5)-PGD & 17 & \textbf{3.20e-06} & 0.06\\
  & MPE(5)-PGD  & 16 & \textbf{7.49e-06}  & 0.05 \\
  & RRE(6)-SGD & 15 & 9.62e-06 &0.12    \\ 
  & MPE(6)-SGD & 16 & 8.24e-06 & 0.20 \\
      \hline 
       
  

 \multicolumn{5}{|c|}{$\alpha=5$} \\
        \hline
         &  GNKS & 25 & 1.59e-04  & 0.08  \\ 
   10 &  GNKS(20) & 20 & 7.66e-04  & 0.05  \\ 
    & RRE(7)-PGD & 20 & \textbf{1.53e-05}  & 0.07 \\
  & MPE(7)-PGD   & 19 &  \textbf{1.76e-05} &0.07  \\
   & RRE(6)-SGD & 14 & 4.39e-05 & 0.13   \\ 
  & MPE(6)-SGD & 14 & 3.90e-05 & 0.13   \\ 
        \hline

     \end{tabular}

\end{table}

\noindent To further validate the efficiency of the proposed hybrid approaches for solving~(\ref{eq2dbratu}), we investigate the effect of the problem dimension. In these experiments, the grid size $n$ is varied from 100 to 1000 in each spatial direction, with several intermediate values considered. For clarity, only representative grid sizes are reported in the figures, since the corresponding curves nearly overlap for closely spaced values of $n$.

\begin{figure}[htbp]
    \centering

    \begin{subfigure}[t]{\textwidth}
        \centering
        \includegraphics[width=3cm,height=5cm]{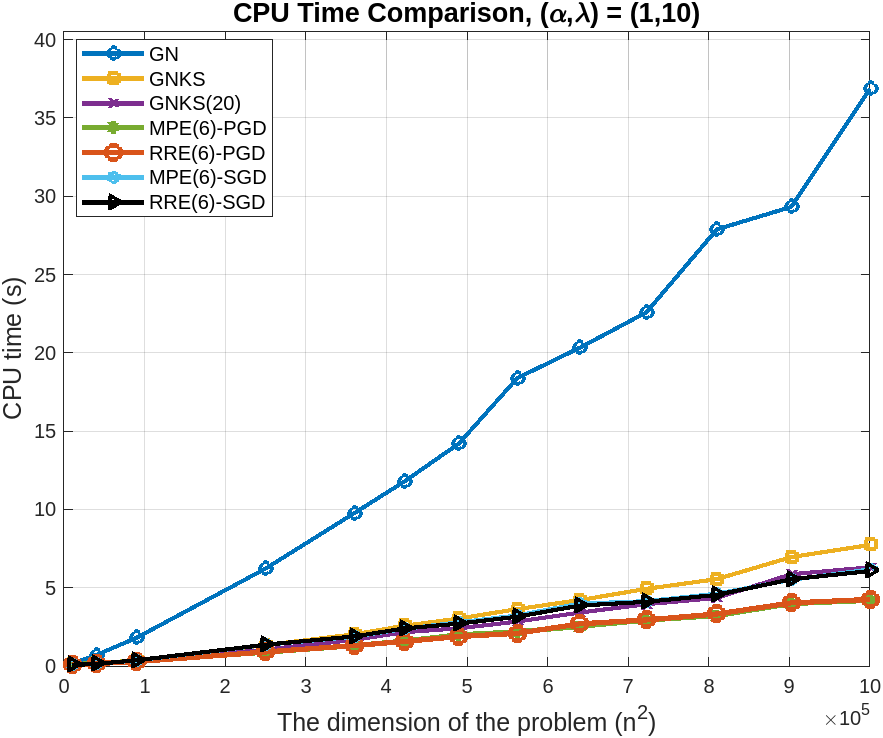}
        \includegraphics[width=3cm,height=5cm]{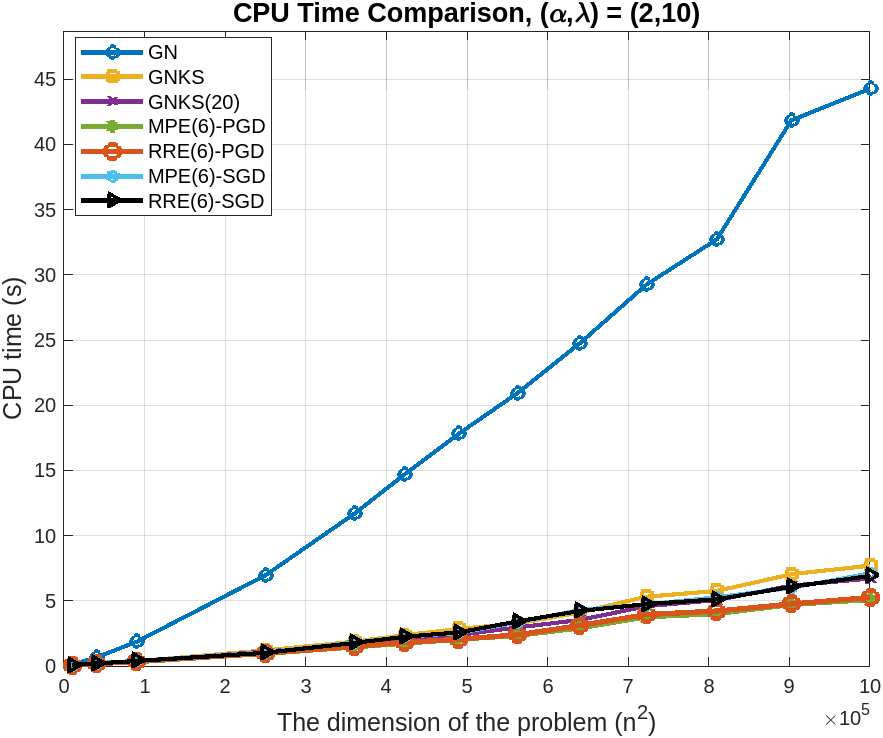}
        \includegraphics[width=3cm,height=5cm]{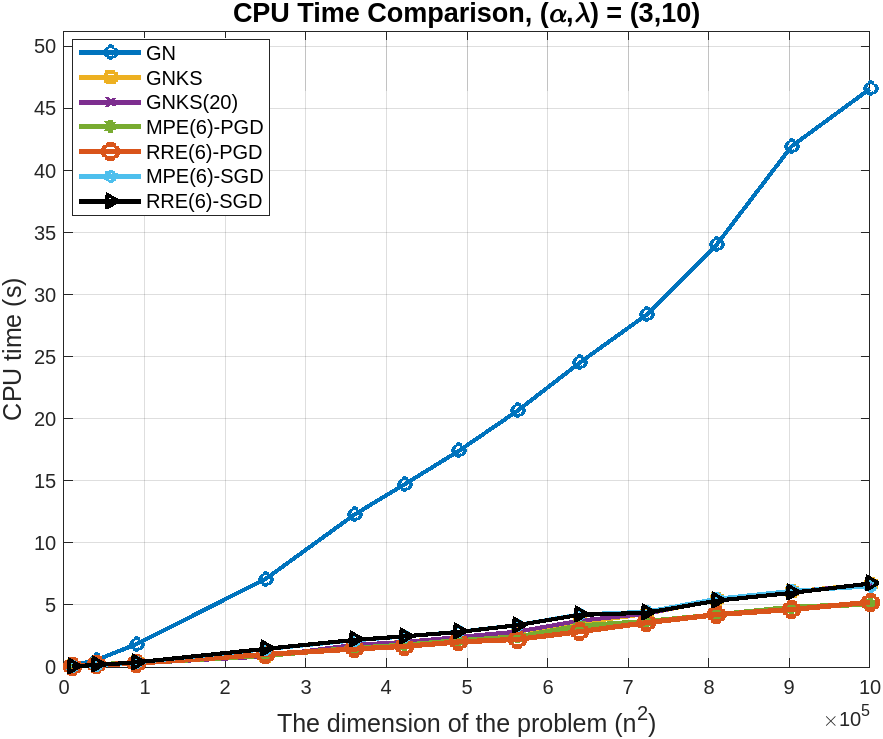}
        \includegraphics[width=3cm,height=5cm]{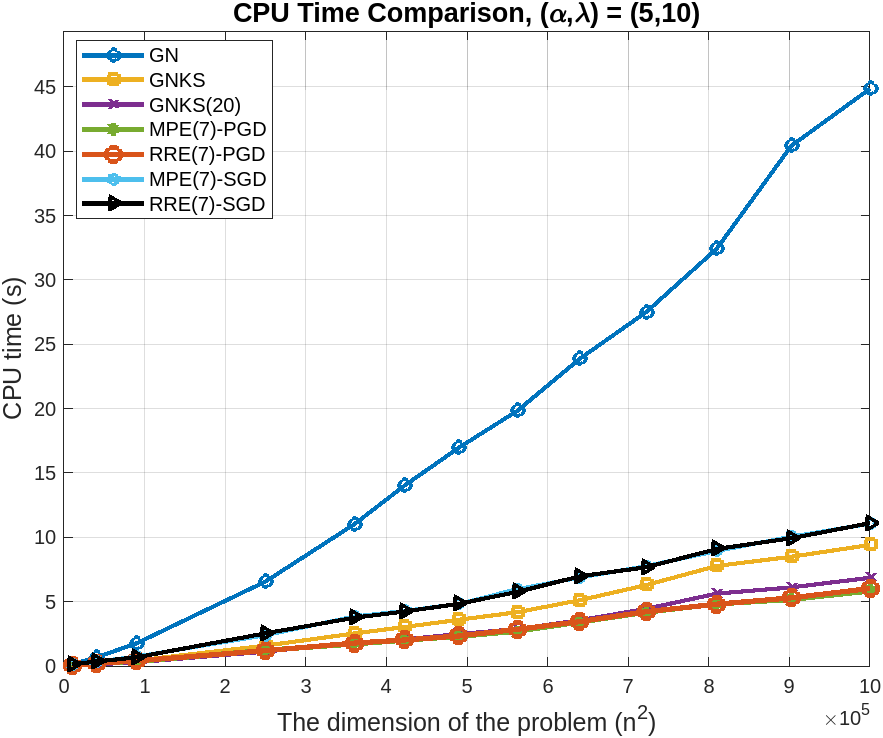}
        \caption{CPU time (in seconds)}
        \label{fig:cpu}
    \end{subfigure}

    \vspace{0.4cm}

    \begin{subfigure}[t]{\textwidth}
        \centering
        \includegraphics[width=3cm,height=5cm]{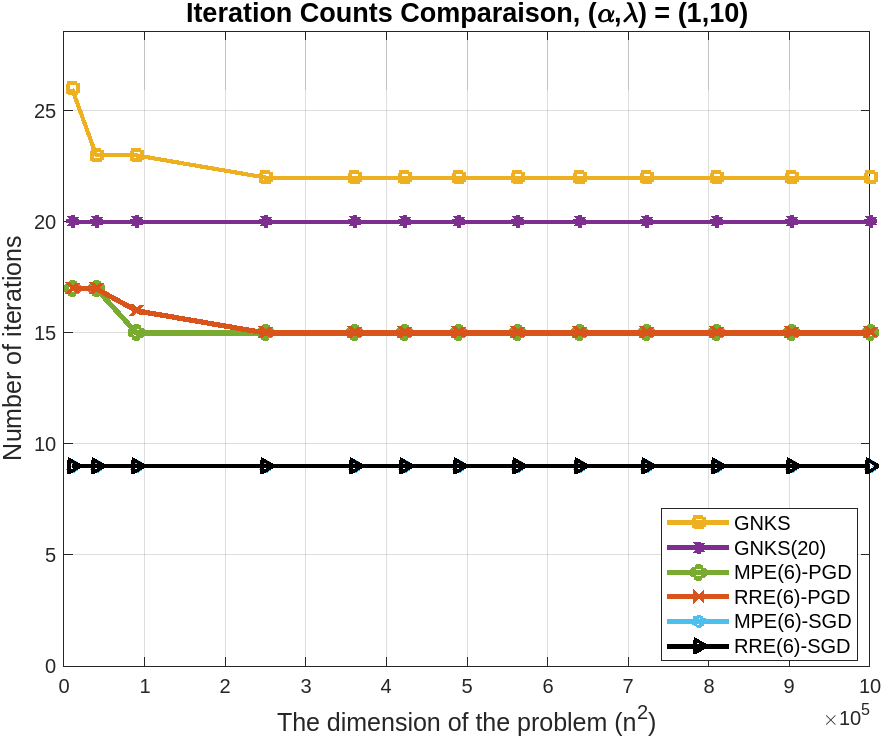}
        \includegraphics[width=3cm,height=5cm]{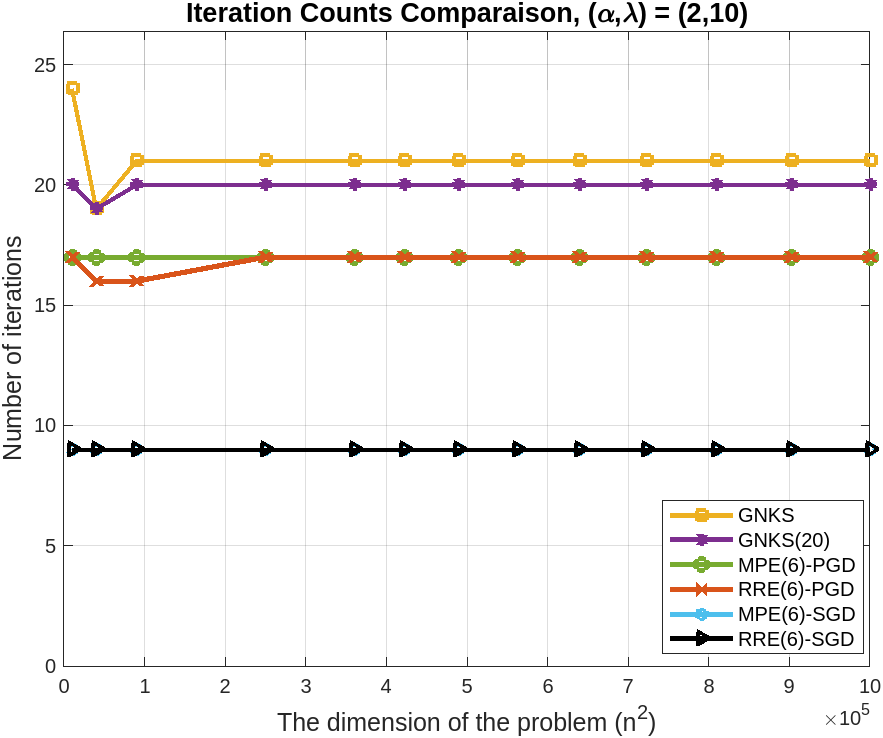}
        \includegraphics[width=3cm,height=5cm]{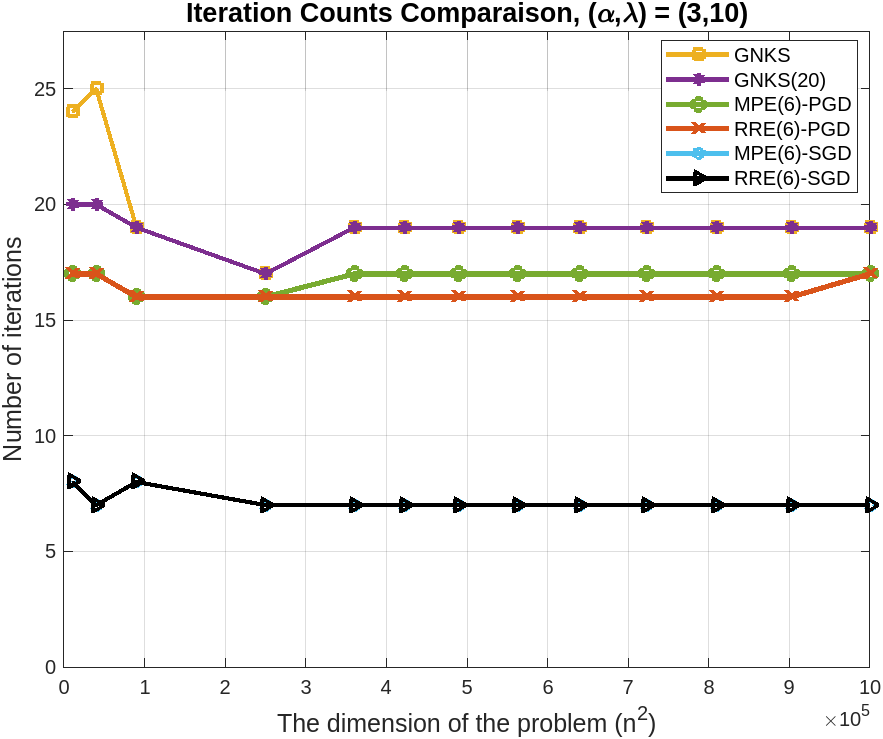}
        \includegraphics[width=3cm,height=5cm]{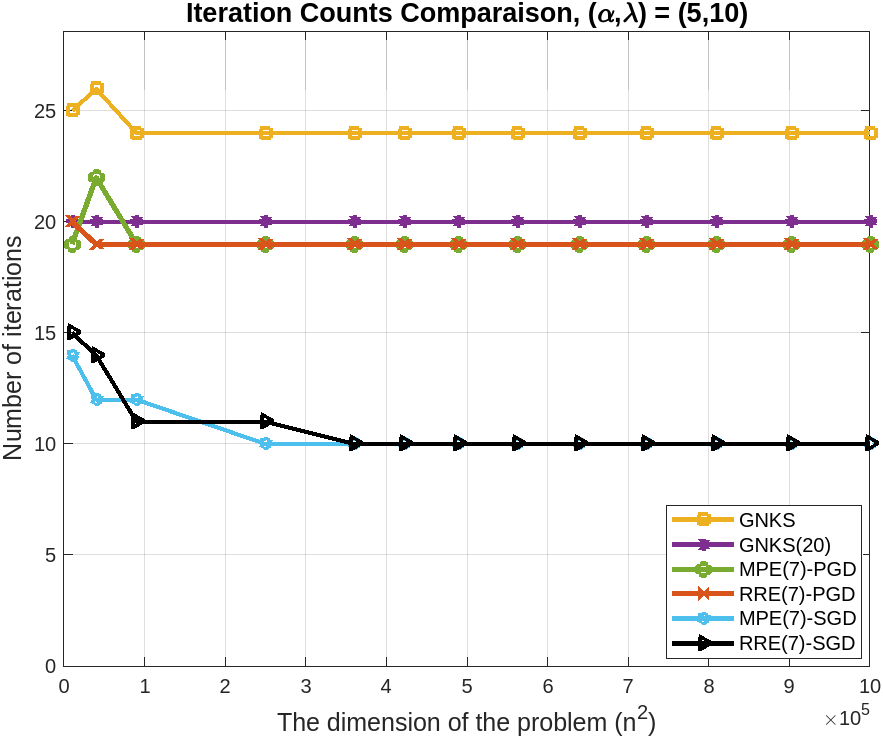}
        \caption{Number of iterations}
        \label{fig:iter}
    \end{subfigure}

    \vspace{0.4cm}

    \begin{subfigure}[t]{\textwidth}
        \centering
        \includegraphics[width=3cm,height=5cm]{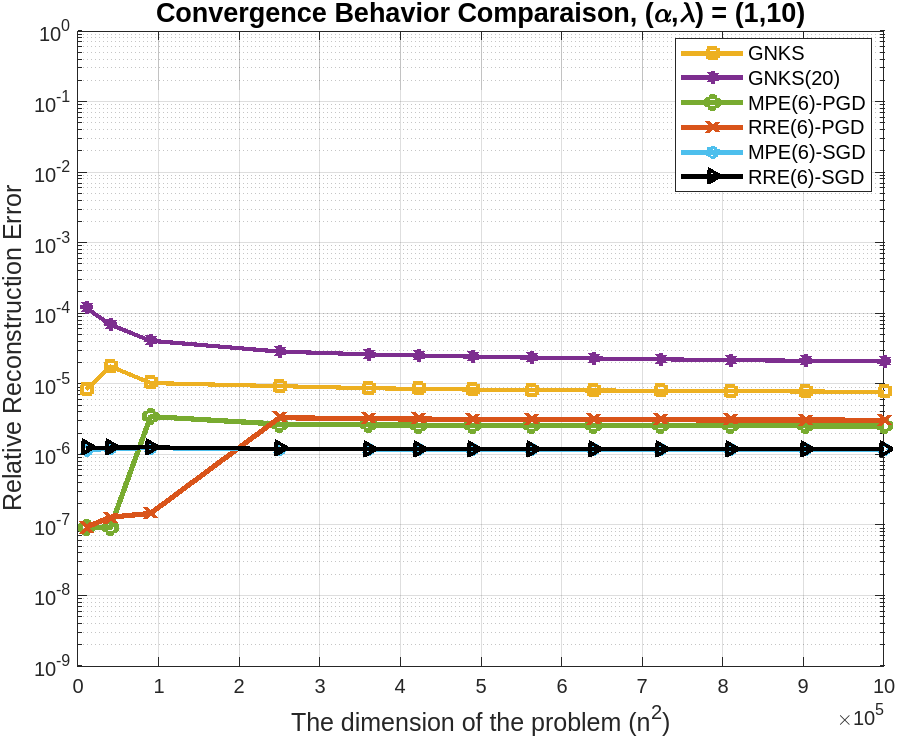}
        \includegraphics[width=3cm,height=5cm]{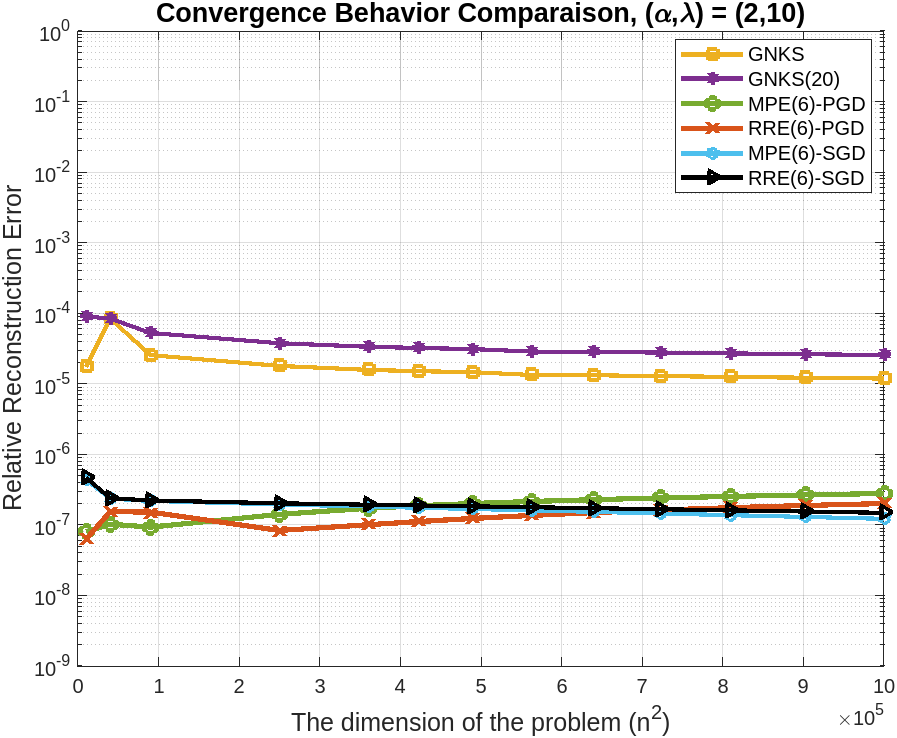}
        \includegraphics[width=3cm,height=5cm]{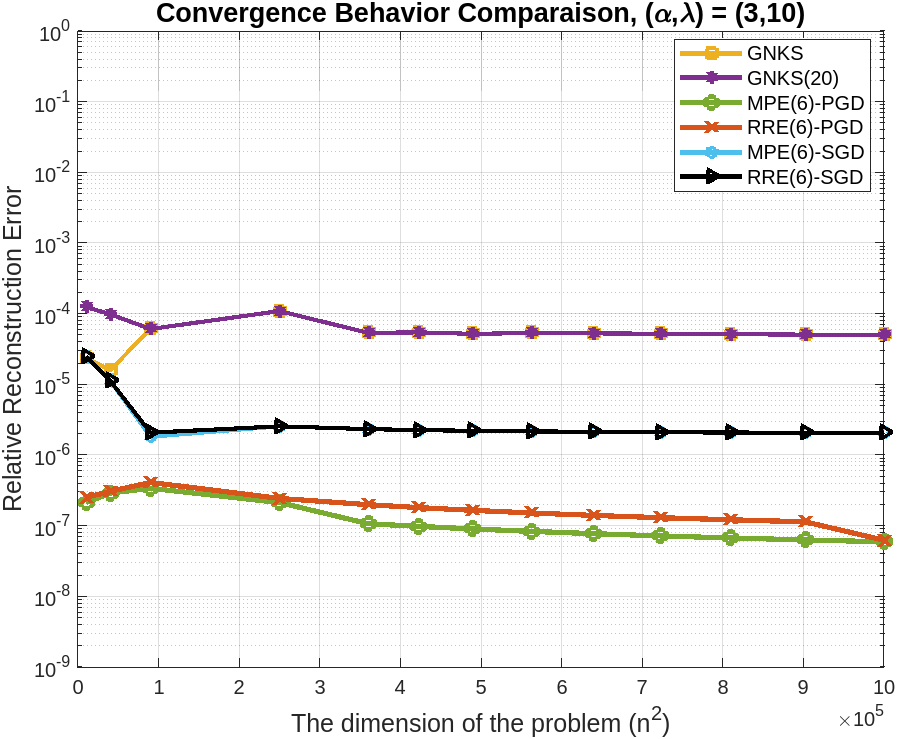}
        \includegraphics[width=3cm,height=5cm]{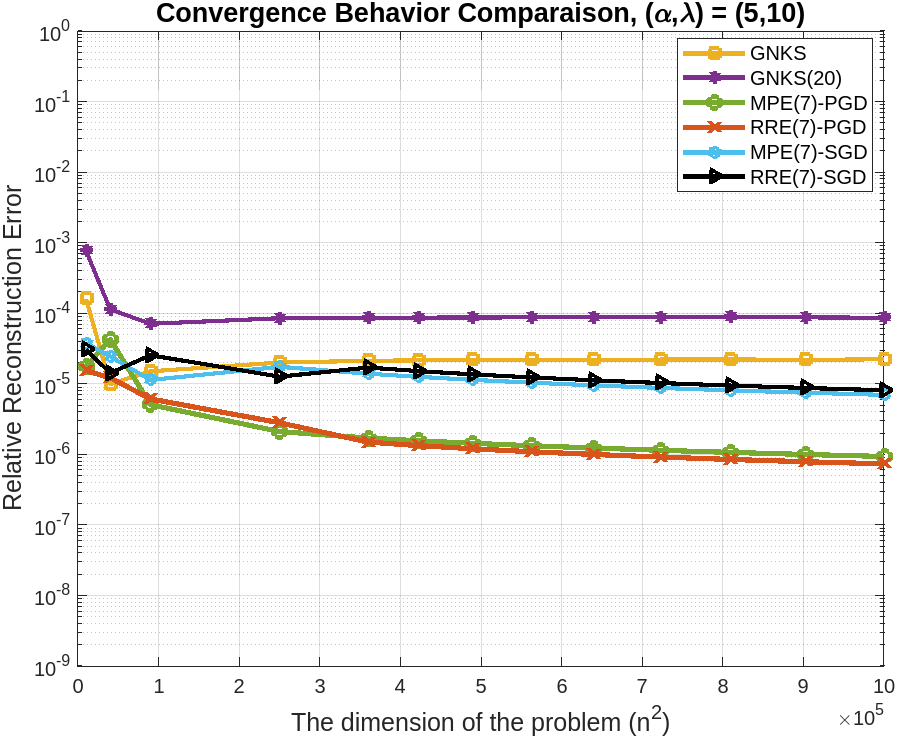}
        \caption{Relative reconstruction error}
        \label{fig:error}
    \end{subfigure}

    \caption{Performance comparison of the different algorithms for increasing dimensions of the 2D Bratu problem with various values of $\alpha$ and $\lambda$.}
    \label{fig:bratu_global}
\end{figure}

For selected values of the parameters $\alpha$ and $\lambda$, we then report the performance of the proposed methods in terms of CPU time (in seconds), iteration count, and relative error (RE) for each fixed grid, as shown in Figures~\ref{fig:cpu}, \ref{fig:iter}, and~\ref{fig:error}, respectively.

\noindent Figure~\ref{fig:cpu} shows that the CPU time required to solve the Bratu problem~(\ref{eq2dbratu}) increases as the grid size $n$ grows. However, when polynomial extrapolation methods (MPE and RRE) are combined with the preconditioned gradient descent approach (PGD), the growth rate of the CPU time is significantly lower than that observed for the other methods, namely GN, GNKS, GNKS(20), MPE-SGD, and RRE-SGD, especially for larger problem dimensions. In addition, the performance of the RRE-PGD and MPE-PGD methods is very similar.

Figures~\ref{fig:iter} and~\ref{fig:error} further indicate that applying restarted polynomial extrapolation techniques (MPE-PGD, RRE-PGD, MPE-SGD, and RRE-SGD) to the sequences generated by the preconditioned gradient descent methods (PGD and SGD) yields more accurate approximate solutions (Figure~\ref{fig:error}) in fewer iterations (Figure~\ref{fig:iter}) than those obtained with the Gauss--Newton method based on generalized Krylov subspaces~\cite{buccini2023efficient} (GNKS and GNKS(20)). Moreover, the MPE-SGD and RRE-SGD methods converge in a substantially smaller number of iterations than the MPE-PGD and RRE-PGD methods. Nevertheless, since MPE-SGD and RRE-SGD require higher computational costs (see Figure~\ref{fig:cpu}), the MPE-PGD and RRE-PGD methods are generally preferable for solving~(\ref{eq2dbratu}).

In the experiments reported in this subsection, we focus on a well-conditioned regime characterized by small values of $\alpha$ and relatively large values of $\lambda$. In this setting, polynomial extrapolation combined with preconditioned gradient-type fixed-point iterations (PGD and SGD) leads to improved convergence behavior compared with both GNKS and its restarted variant GNKS($q=20$). In particular, the proposed hybrid approaches (RRE-PGD, MPE-PGD, RRE-SGD, and MPE-SGD) achieve a significant reduction in the number of iterations while producing more accurate reconstructions, as measured by the relative error (RE).

In contrast, in an ill-conditioned regime characterized by larger values of $\alpha$ and smaller values of $\lambda$, polynomial extrapolation applied to GD, PGD, and SGD iterations does not guarantee a systematic reduction in the iteration count, nor does it consistently improve the reconstruction accuracy.

In this regime, the restarted GNKS method with restart parameter $q=20$ proves to be a robust and efficient strategy. Restarting significantly reduces the number of iterations required to reach a prescribed tolerance on the relative residual norm, while preserving essentially the same reconstruction accuracy as the non-restarted GNKS method. Consequently, GNKS($q=20$) provides a reliable and stable baseline in ill-conditioned settings.

As an alternative, polynomial extrapolation can be applied to the sequence of GNKS iterates. Our numerical results indicate that this hybrid strategy may further reduce the iteration count compared with GNKS($q=20$). However, this reduction does not translate into an improvement in reconstruction quality, and the resulting relative reconstruction error is generally larger than that obtained with the restarted GNKS approach.

Therefore, in ill-conditioned regimes, we favor the restarted GNKS strategy, which offers a more reliable balance between iteration count, numerical stability, and reconstruction accuracy.

\subsubsection{The standard Bratu problem}
We now perform a numerical study to assess the effectiveness of the proposed techniques for the standard two-dimensional Bratu problem, corresponding to $\alpha = 0$ in~(\ref{eq2dbratu}). In this setting, we focus on challenging regimes characterized by large values of the parameter $\lambda$, which are well known to induce strong nonlinear effects and increased numerical difficulty.

Such parameter choices lead to a more demanding computational regime, often associated with slow convergence and reduced robustness of classical iterative solvers. Within this context, Table~\ref{table_initila_bratu} reports the performance of the proposed methods and illustrates their ability to improve the accuracy of the computed approximations for the Bratu problem.

\begin{table}[htbp]
	 \caption{Performance of the proposed methods for the standard Bratu problem $-\Delta x +\lambda e^x=y$ for high values of $\lambda$.}
  \label{table_initila_bratu}	
  \centering
    \begin{tabular}{|c|c|c|c|c|}
        \hline
  & Method  & Iter & Relative Error (RE) &  CPU(s)  \\        
        \hline
     &  GNKS & 19 & 3.01e-05  & 0.08 \\ 
   &  GNKS(10) & 10 &  2.39e-03 & 0.04 \\  
 $\lambda=10^{1}$  & VEA(5)-PGD   & 19 &  1.62e-06 & 0.07 \\  
    & RRE(5)-PGD & 15 & 2.71e-06 &  0.06\\ 
  & MPE(5)-PGD  & 15 & 2.48e-06 & 0.07 \\   
   & VEA(5)-SGD & 12 & 9.80e-07 & 0.10 \\ 
   & RRE(5)-SGD & 10 & 6.33e-06 &  0.08 \\ 
  & MPE(5)-SGD  & 10 & 7.80e-06 & 0.09 \\ 
  
\hline
   &  GNKS & 13 & 4.91e-06 & 0.04 \\ 
   &  GNKS(10) & 10 & 7.82e-05 & 0.02 \\  
 $\lambda= 10^{4}$  & VEA(5)-PGD   & 12 & 3.13e-07  & 0.10 \\  
    & RRE(5)-PGD & 14 & 4.46e-08 &  0.12\\ 
  & MPE(5)-PGD  & 14 & 5.47e-08 & 0.13  \\ 
  
   & VEA(2)-SGD  & 8 &  \textbf{3.51e-10} & 0.09 \\ 
   & RRE(2)-SGD & 7 & \textbf{3.08e-12} & 0.08 \\ 
  & MPE(2)-SGD   & 6 &  \textbf{4.94e-09} &  0.07\\ 
\hline
   &  GNKS & 11 &  5.55e-06 & 0.03 \\ 
   &  GNKS(10) & 10 &  7.79e-05 &  0.02\\  
$\lambda= 10^{5}$   & VEA(5)-PGD   & 23 & 7.12e-09 & 0.19 \\  
    & RRE(5)-PGD & 15 & 4.22e-08 &  0.14\\ 
  & MPE(5)-PGD  & 15 & 4.45e-08 & 0.17 \\ 
  
    & VEA(2)-SGD &  8 & \textbf{1.46e-13}  & 0.09 \\ 
   & RRE(2)-SGD & 7 &  \textbf{3.75e-14} & 0.07 \\ 
  & MPE(2)-SGD  & 7 & \textbf{3.39e-14} & 0.08 \\ 
   
\hline
   &  GNKS & 11 & 1.56e-06 & 0.03 \\ 
  &  GNKS(10) & 10 & 7.80e-05 & 0.02 \\  
$\lambda= 10^{6}$   & VEA(5)-PGD   &  12 & 1.25e-06 &  0.15\\  
    & RRE(5)-PGD & 13 & 1.49e-06 &   0.18\\ 
  & MPE(5)-PGD  & 13 & 1.44e-06 & 0.14  \\ 
  
   & VEA(2)-SGD  & 9 & \textbf{1.04e-15} & 0.10 \\ 
   & RRE(2)-SGD &  7 & \textbf{1.15e-15} & 0.08\\ 
  & MPE(2)-SGD   & 7 & \textbf{1.04e-15} &  0.09 \\
   
\hline

\end{tabular}
  
\end{table}

According to Table~\ref{table_initila_bratu}, the restarted GNKS(20) method improves the convergence of the GNKS method without restarting in terms of the number of iterations and CPU time. However, it yields less accurate approximate solutions for the standard Bratu problem than GNKS. More precisely, the $L_2$ norms of the relative reconstruction error (RE) increase when GNKS(20) is employed.

\begin{figure}[h!]
    \centering

    \begin{subfigure}[t]{\textwidth}
        \centering
        \includegraphics[width=4.4cm,height=5cm]{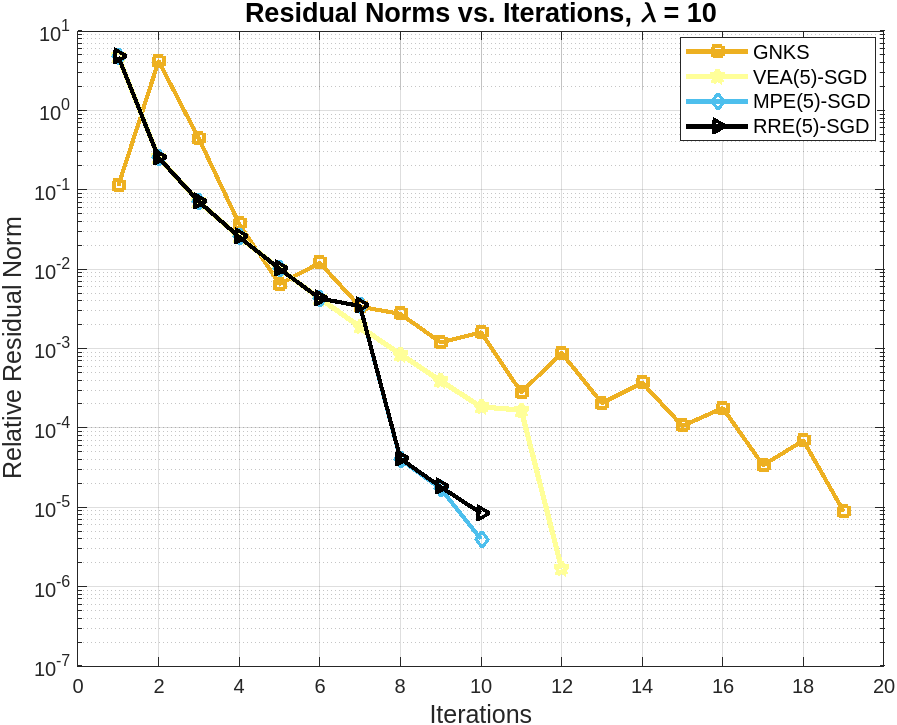}
        \includegraphics[width=4.4cm,height=5cm]{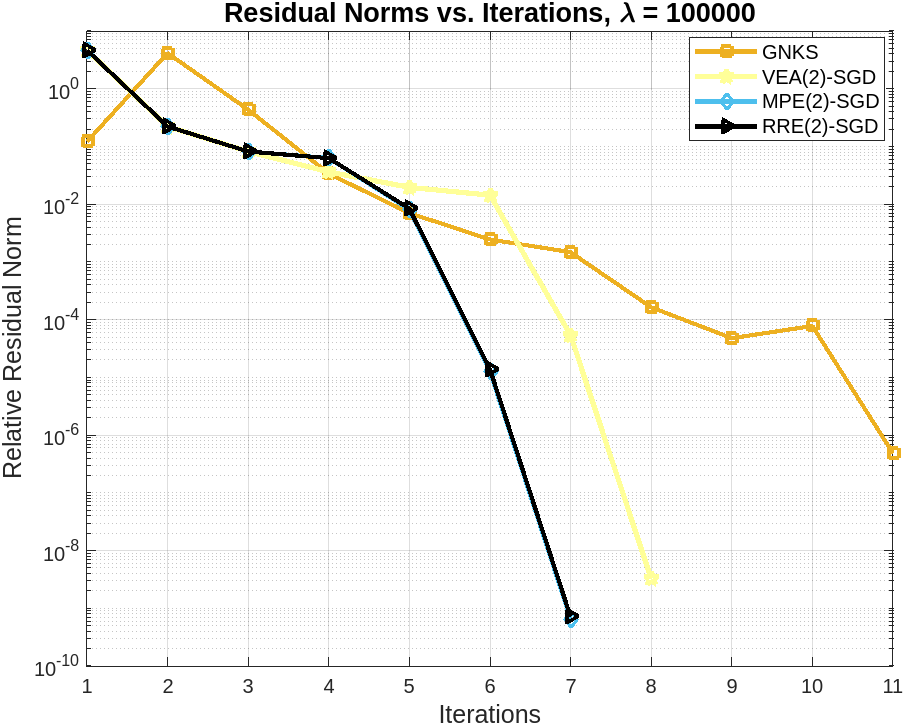}
        \includegraphics[width=4.4cm,height=5cm]{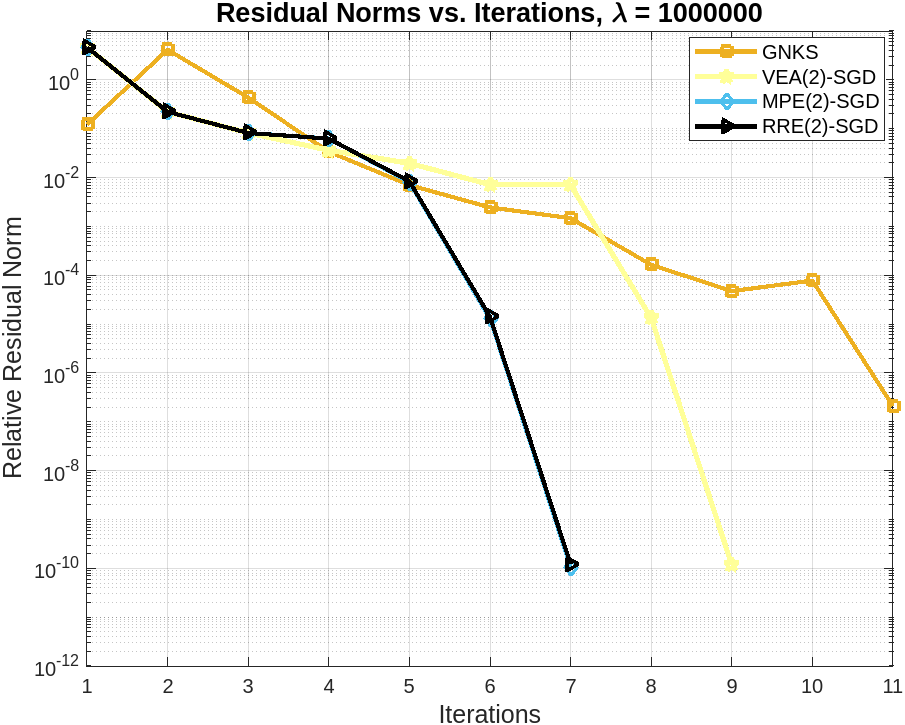}
        \caption{Relative residual norms}
        \label{fig:bratu_residual}
    \end{subfigure}

    \vspace{0.4cm}

    \begin{subfigure}[t]{\textwidth}
        \centering
        \includegraphics[width=4.4cm,height=5cm]{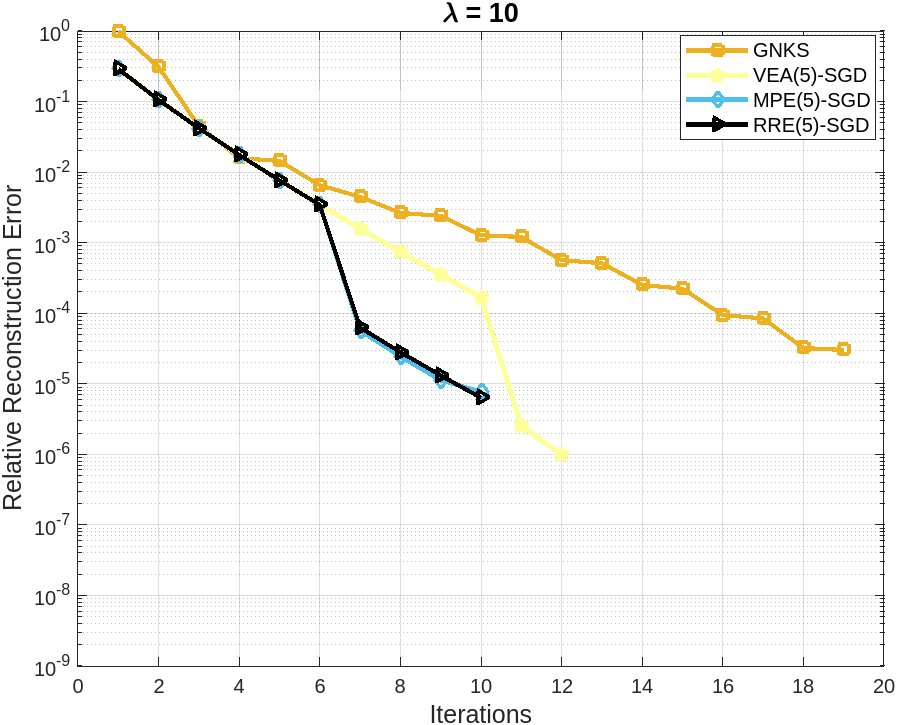}
        \includegraphics[width=4.4cm,height=5cm]{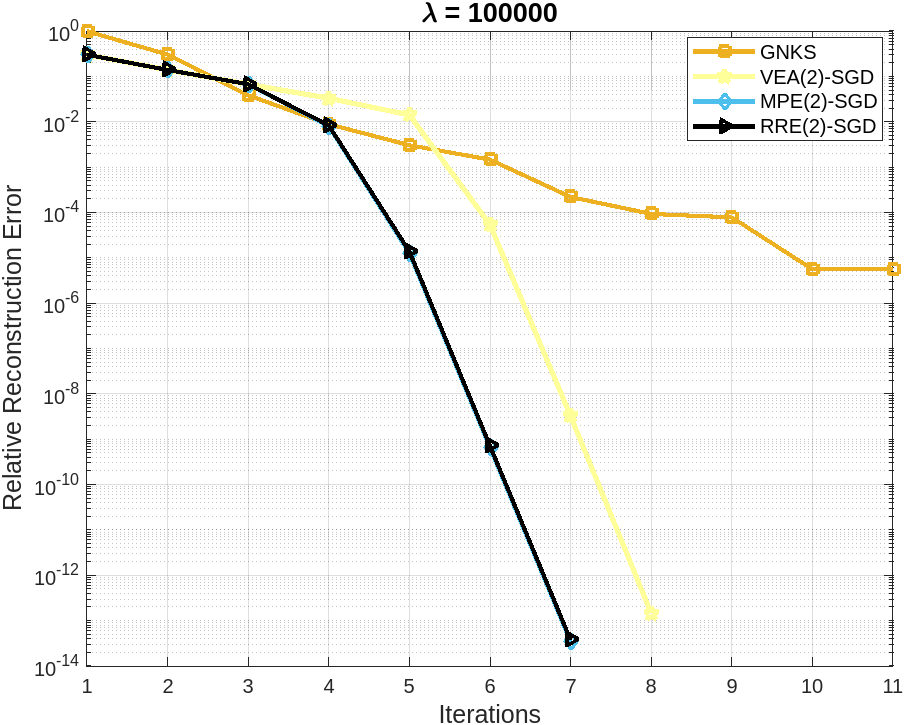}
        \includegraphics[width=4.4cm,height=5cm]{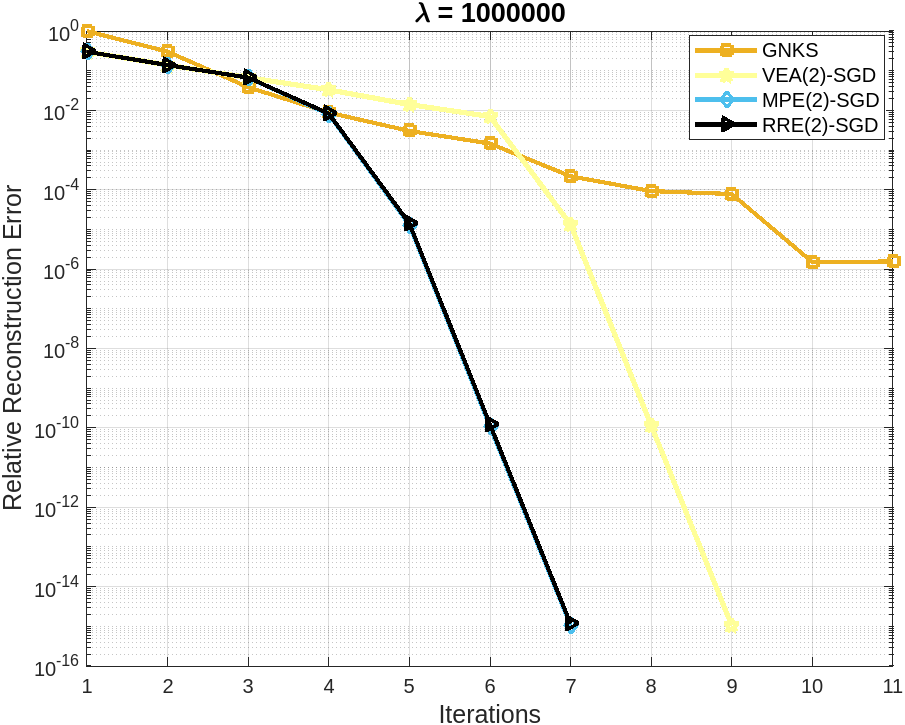}
        \caption{Relative reconstruction errors}
        \label{fig:bratu_error}
    \end{subfigure}

    \caption{Behavior of the VEA, RRE, and MPE methods applied to SGD iterations for the standard 2D Bratu problem
    $-\Delta x + \lambda e^{x} = y$ with high values of $\lambda$.}
    \label{fig:BratuHighLambda}
\end{figure}

By contrast, our hybrid approaches, based on the application of the vector $\epsilon$-algorithm and polynomial extrapolation methods to the PGD and SGD iterations, achieve improved convergence behavior compared with both GNKS and GNKS(20). In particular, a noticeable reduction in the relative error (Table~\ref{table_initila_bratu} and Figure~\ref{fig:BratuHighLambda}) is obtained after only a small number of iterations. This improvement is especially pronounced for large values of the parameter $\lambda$ in the standard Bratu problem, where the VEA-SGD, RRE-SGD, and MPE-SGD methods outperform their PGD-based methods, namely VEA-PGD, RRE-PGD, and MPE-PGD.

Consequently, the standard Bratu problem can be solved more effectively with the proposed hybrid approaches than with the Gauss--Newton method based on generalized Krylov subspaces~\cite{buccini2023efficient} (GNKS), yielding more accurate approximate solutions, particularly for very large values of $\lambda$ corresponding to highly challenging regimes, as illustrated in Table~\ref{table_initila_bratu} and Figure~\ref{fig:BratuHighLambda}.

\subsection{An extremely sparse nonlinear problem}

We now consider an extremely sparse problem (see~\cite{buccini2023efficient}) defined by a differentiable nonlinear function 
$f : \mathbb{R}^{n} \rightarrow \mathbb{R}^{n-1}$, whose $i$-th component is given by
\[
f_i(\mathrm{x}) = \sin(x_i + x_{i+1}), \qquad i = 1, 2, \ldots, n-1.
\]

\noindent The associated bidiagonal Jacobian $J_f(\mathrm{x}) \in \mathbb{R}^{(n-1)\times n}$ is given by
\[
J_f(\mathrm{x}) =
\begin{bmatrix}
\cos(x_1 + x_2) & \cos(x_1 + x_2) & & & \\
& \cos(x_2 + x_3) & \cos(x_2 + x_3) & & \\
& & \ddots & \ddots & \\
& & & \cos(x_{n-1} + x_n) & \cos(x_{n-1} + x_n)
\end{bmatrix}.
\]

\noindent The sparse problem is generated by selecting the exact solution $\mathrm{x} = \tfrac{1}{2}\sin(x), \\\quad x \in (-\pi,\pi),$
sampled on an equispaced grid of size 
$n \in \{10^{3}, 10^{4}, 10^{5}, 10^{6}, 10^{7}\}$, which defines the vector $\mathrm{x} \in \mathbb{R}^n$ of the minimization problem~(\ref{Lstsq}).

Table~\ref{table_sparse_prob} reports the convergence results obtained with the Gauss--Newton (GN) method, the Gauss--Newton method based on generalized Krylov subspaces~\cite{buccini2023efficient} (GNKS), its restarted variant GNKS($q$) with restart parameter $q=10$, and the restarted polynomial extrapolation methods MPE($q$)-SGD and RRE($q$)-SGD applied to the sequence generated by the scaled gradient descent method for different grid sizes $n$. For each fixed grid size, we report the number of iterations (“Iter”) required to reduce the relative residual norm $\frac{\Vert \mathrm{x}_{k+1} - \mathrm{x}_k \Vert_2}{\Vert \mathrm{x}_k \Vert_2}$
below $10^{-5}$. The relative reconstruction error (“RE”),$
\frac{\Vert \mathrm{x} - \mathrm{x}_{\text{true}} \Vert_2}{\Vert \mathrm{x}_{\text{true}} \Vert_2},$
as well as the CPU time in seconds (“CPU(s)”), are also reported in Table~\ref{table_sparse_prob}.

\begin{table}[htbp]
  \caption{Performance of GN, GNKS, GNKS($q$), VEA($q$)-SGD, MPE($q$)-SGD, and RRE($q$)-SGD for solving the extremely sparse problem for different grids.}
  \label{table_sparse_prob}
  \centering
 \begin{tabular}{|c|c|c|c|}
        \hline
       Method  & Relative Error (RE) &  Iter  &  CPU(s)  \\
        \hline
         \multicolumn{4}{|c|}{$n=10^{3}$} \\
        \hline
            GN & 8.50e-02 & 5  & 0.00107 \\ 
         GNKS & 1.22e-04 & 17 & 0.0145 \\  
       GNKS(10) & 1.50e-03 & 10 &  0.0137\\
   RRE(1)-SGD &  6.68e-05 &  9  & 0.0165 \\ 
  MPE(1)-SGD   & 6.68e-05 & 9  & 0.0159 \\ 
   VEA(1)-SGD &  6.31e-05  & 15 & 0.0212 \\ 
  RRE(3)-SGD & 5.88e-05  & 9 & 0.0139  \\ 
  MPE(3)-SGD   & 6.38e-05 & 8  & 0.0123 \\ 
VEA(3)-SGD & 5.53e-05  & 16 & 0.0220 \\ 
  RRE(5)-SGD &  4.92e-05 & 11 & 0.0202 \\ 
  MPE(5)-SGD   & 4.92e-05 & 11  & 0.0169  \\ 
 VEA(5)-SGD & 5.16e-05  & 14& 0.0272 \\



         
         

  
       \hline 
        \multicolumn{4}{|c|}{$n=10^{6}$} \\
        \hline
        GN & 2.82e-03  & 5  & 0.620  \\ 
         GNKS  &  7.26e-06 & 9  & 1.036 \\ 
         
   RRE(1)-SGD  & 2.22e-09  &  10 & 1.796  \\ 
 MPE(1)-SGD    & 2.22e-09  &  10 & 2.031 \\ 
VEA(1)-SGD &  1.88e-09 & 31 &  6.751s\\ 
  RRE(5)-SGD & 1.91e-09  & 12 & 3.237 \\ 
  MPE(5)-SGD   & 2.01e-09 & 12  & 3.111 \\ 
VEA(5)-SGD &  1.85e-09 & 23 &  5.864\\ 
  RRE(6)-SGD & 2.27e-09  & 9 & 2.572 \\ 
  MPE(6)-SGD  & 2.27e-09 & 9  & 2.520 \\ 
VEA(6)-SGD &  2.18e-09 & 26 & 6.459\\ 

       \hline 
 \multicolumn{4}{|c|}{$n=10^{7}$} \\
        \hline
         GN & 8.94e-04  & 5  & 11.707 \\ 
         GNKS  &  2.01e-06 &  8 & 16.398 \\ 
           RRE(1)-SGD  & 7.02e-11  & 10  & 31.530 \\ 
 MPE(1)-SGD    & 7.02e-11  &   10 & 31.566 \\ 
VEA(1)-SGD &  7.25e-06 & 35 &  106.801\\ 
  
RRE(7)-SGD & 6.96e-11  & 10 & 45.770 \\ 
  MPE(7)-SGD   &   6.96e-11 & 10   &  43.472 \\ 
VEA(4)-SGD &  5.72e-11 & 28 &  115.419\\

       \hline 
  
 \end{tabular}
\end{table}
 
\noindent Table~\ref{table_sparse_prob} shows that, compared with the GN method, the GNKS approach yields significantly smaller relative reconstruction errors (RE) for all grid sizes. This indicates that GNKS outperforms GN in terms of approximation accuracy, although the GN method converges with substantially lower CPU time and fewer iterations than GNKS. Despite its higher computational cost, GNKS remains preferable when accuracy is the primary concern. In addition, the restarted version of GNKS with $q=10$ produces less accurate solutions in terms of relative error compared with the non-restarted GNKS method.

On the other hand, combining the polynomial extrapolation methods MPE and RRE with the SGD iteration results in more accurate approximate solutions than GNKS, while requiring fewer iterations and maintaining acceptable CPU times. This improvement becomes more pronounced as the grid size increases. Consequently, we favor the MPE-SGD and RRE-SGD approaches for solving this extremely sparse problem.

\begin{figure}[h!]
    \centering

    \begin{subfigure}[t]{0.48\textwidth}
        \centering
        \includegraphics[width=7cm,height=7cm]{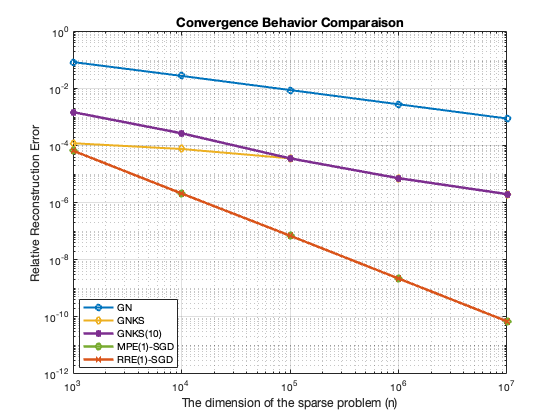}
        \caption{Error norms for different grid sizes}
        \label{fig:sparse_grids}
    \end{subfigure}
    \hfill
    \begin{subfigure}[t]{0.48\textwidth}
        \centering
        \includegraphics[width=7cm,height=7cm]{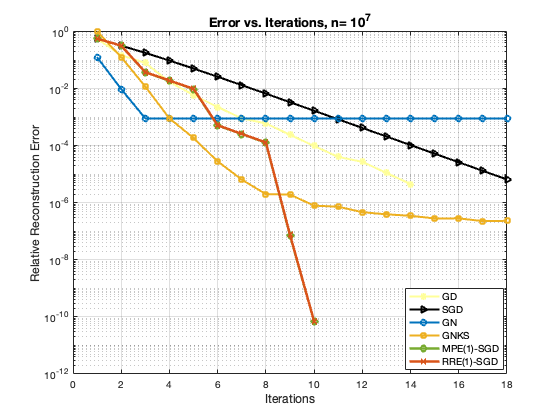}
        \caption{Fixed grid size, $n = 10^{7}$}
        \label{fig:sparse_fixed}
    \end{subfigure}

    \caption{Convergence comparison between the GN and GNKS methods  with SGD using acceleration techniques for an extremely sparse problem.}
    \label{fig:sparse_gn_gnks}
\end{figure}

Figure~\ref{fig:sparse_gn_gnks} provides a graphical illustration of the results reported in Table~\ref{table_sparse_prob}. In Figure~\ref{fig:sparse_grids}, we examine the convergence behavior of the different methods, including GN, GNKS, GNKS(10), MPE(1)-SGD, and RRE(1)-SGD, by plotting the relative reconstruction error as a function of the problem dimension. The results indicate that the MPE-SGD and RRE-SGD methods consistently outperform the other approaches in terms of solution accuracy for all grid sizes.

Furthermore, Figure~\ref{fig:sparse_fixed} shows the relative reconstruction error as a function of the iteration number for a high-dimensional sparse problem with $n=10^{7}$. In this case, we also include gradient descent (GD) and scaled gradient descent (SGD) for comparison, in addition to GN and GNKS. These results demonstrate that polynomial extrapolation significantly enhances the convergence of gradient-based methods and confirm the efficiency of the proposed approaches compared with all other methods, particularly in very high-dimensional settings.

\section{Conclusion}\label{sec5}
To demonstrate the effectiveness of vector extrapolation methods for solving nonlinear problems, this work focused on nonlinear least-squares problems by combining preconditioned gradient-based methods (PGD and SGD) with vector extrapolation techniques, including polynomial-type extrapolation methods (RRE and MPE) and the vector $\varepsilon$-algorithm (VEA).

A comprehensive numerical study was first conducted on the Bratu problem in its general form for various values of the parameters $\alpha$ and $\lambda$. The results show that applying restarted (cyclic) polynomial extrapolation methods and the vector $\varepsilon$-algorithm to preconditioned gradient iterations leads to superior convergence behavior in terms of relative reconstruction error when compared with Gauss–Newton methods based on generalized Krylov subspaces (GNKS) and their restarted variants. In particular, although restarting GNKS reduces the number of iterations required for convergence, it does not improve the accuracy of the approximate solution relative to the non-restarted GNKS method.

Since restarted extrapolation strategies naturally account for memory and storage constraints, polynomial extrapolation methods were preferred within this framework. Their effectiveness was further validated on the standard Bratu problem in a challenging regime corresponding to large values of the parameter $\lambda$. In this setting, the proposed approaches consistently produced more accurate approximate solutions, especially for higher values of $\lambda$, while requiring fewer iterations than Gauss–Newton methods based on generalized Krylov subspaces.

Finally, an extremely sparse nonlinear problem was considered. For this problem, GNKS yielded more accurate solutions than the standard Gauss–Newton (GN) method, while GN converged with significantly lower CPU time and fewer iterations. In contrast, the proposed extrapolation-accelerated solvers achieved higher accuracy than GNKS with comparable computational costs and reduced iteration counts. Moreover, the improvement in convergence, measured in terms of relative reconstruction error, became increasingly pronounced as the problem dimension increased.

Overall, these results demonstrate that combining vector extrapolation techniques with preconditioned gradient-based methods provides a competitive, accurate, and robust alternative for the solution of large-scale nonlinear least-squares problems.

\bibliographystyle{elsarticle-num-names}
\bibliography{bibliography}



\end{document}